\documentclass[11pt]{article}

\usepackage{color}
\usepackage{latexsym}
\usepackage{amssymb}
\usepackage{amsmath, amsfonts,amssymb,theorem,euscript,array,enumerate,amsfonts,mathrsfs}
\usepackage{graphicx}

\newtheorem{Theorem}{Theorem}[part]

\newtheorem{Proposition}{Proposition}[part]

\newtheorem{Lemma}{Lemma}[part]
\newtheorem{Corollary}{Corollary}[part]
\newtheorem{Remark}{Remark}[part]

\def\esssup_#1{\underset{#1}{\mathrm{ess\,sup\, }}}
\def\essinf_#1{\underset{#1}{\mathrm{ess\,inf\, }}}

\def \trans{^{\scriptscriptstyle{\intercal}}}

\def \Frac{\displaystyle\frac}

\def \trans{^{\scriptscriptstyle{\intercal }}}

\def \N{\mathbb{N}}
\def \R{\mathbb{R}}

\def \E{\mathbb{E}}
\def \F{\mathbb{F}}

\def \P{\mathbb{P}}

\def \Ac{{\cal A}}
\def \Bc{{\cal B}}

\def \Ec{{\cal E}}
\def \Fc{{\cal F}}

\def \Lc{{\cal L}}
\def \Pc{{\cal P}}

\def \Oc{{\cal O}}

\def \Tc{{\cal T}}

\def \Vc{{\cal V}}

\def \Vc{{\cal V}}

\def \ni{\noindent}

\def \eps{\varepsilon}
\def \ups{\upsilon}
\def \var{\vartheta}

\def \ep{\hbox{ }\hfill$\Box$}

\def\Dt#1{\Frac{\partial #1}{\partial t}}

\def\reff#1{{\rm(\ref{#1})}}

\def\beqs{\begin{eqnarray*}}
\def\enqs{\end{eqnarray*}}
\def\beq{\begin{eqnarray}}
\def\enq{\end{eqnarray}}

\allowdisplaybreaks

\begin{document}

\title{Reflected BSDEs with nonpositive jumps, and controller-and-stopper games\thanks{We would like to thank  Said Hamad\`ene  for helpful conversations during the preparation of this work.}}

\author{S\'ebastien CHOUKROUN\thanks{Laboratoire de Probabilit\'es et Mod\`eles Al\'eatoires, CNRS, UMR 7599, Universit\'e Paris Diderot,
		\sf  choukroun at math.univ-paris-diderot.fr}~~~
		Andrea COSSO\thanks{Dipartimento di Matematica, Politecnico di Milano,
               \sf  andrea.cosso at polimi.it }~~~
               Huy{\^e}n PHAM\thanks{Laboratoire de Probabilit\'es et Mod\`eles Al\'eatoires, CNRS, UMR 7599, Universit{\'e} Paris Diderot, and
               CREST-ENSAE,  \sf pham at math.univ-paris-diderot.fr}
             }

\maketitle

\date{}

\begin{abstract}
We study a class of  reflected backward stochastic differential equations with nonpositive jumps and upper barrier. Existence and uniqueness of a minimal solution is proved by a double pena\-lization approach under regularity  assumptions on the obstacle.  In a suitable regime switching diffusion framework, we show the  connection between our class of BSDEs and  fully nonlinear variational inequalities.  Our BSDE representation provides in particular a Feynman-Kac type formula for PDEs associated to  general zero-sum stochastic differential controller-and-stopper games, where control affect both drift and diffusion term, and the diffusion coefficient can be degenerate. Moreover, we state a dual game formula of this BSDE minimal solution involving equivalent change of probability measures, and discount processes. This gives in particular a new representation for zero-sum stochastic differential controller-and-stopper games.
\end{abstract}

\vspace{5mm}

\noindent {\bf Key words:}  Backward stochastic differential equations (BSDE) with constrained jumps, reflected BSDE, regime-switching jump-diffusion, Hamilton-Jacobi-Bellman Isaacs equation,
controller-and-stopper game.

\vspace{5mm}

\noindent {\bf MSC Classification:}  60H10,  60H30, 35K86, 49N70.

\newpage

\section{Introduction}

Backward stochastic differential equations (BSDEs), introduced in the seminal paper by Pardoux and Peng \cite{parpen90},  have emerged over the last years as a major topic in probability, especially through its deep connection with nonlinear PDEs and associated probabilistic numerical methods,  and stochastic control in mathematical finance.  A solution to a standard BSDE on a filtered probability space $(\Omega,\Fc,(\Fc_t)_{0\leq t\leq T},\P)$ generated by an $\R^d$-valued Brownian motion $W$, is a pair of a progressively measurable process $(Y,Z)$ satisfying:
\beq \label{BSDEstandard}
Y_t &=& \xi + \int_t^T F(s,Y_s,Z_s) ds - \int_t^T Z_s dW_s, \;\;\; 0 \leq t\leq T,
\enq
where the generator $F$ is a progressively measurable function, and the terminal data $\xi$ is $\Fc_T$-measurable. In the Markovian case where
$\xi(\omega)$ $=$ $g(W_T(\omega))$,  $F(t,\omega,y,z)$ $=$ $f^0(W_t(\omega),y,z)$, for some continuous functions  $g$ and  $f^0$ on  $\R^d$ and
$\R^d\times\R\times\R^d$,  it is well-known from \cite{parpen92} that BSDE \reff{BSDEstandard} provides a Feynman-Kac formula to the semi-linear partial differential equation (PDE):
\beq \label{semiPDE}
\Dt{v} + \frac{1}{2}{\rm tr}(D_x^2 v) + f^0(x,v,D_x v) &=& 0, \;\;\; \mbox{ on } [0,T)\times\R^d,
\enq
with terminal condition $v(T,\cdot)$ $=$ $g$, through the relation: $Y_t$ $=$ $v(t,W_t)$, $0\leq t\leq T$.  We also notice that when the function $f^0$ is  in the form:  $f^0(x,z)$ $=$ $\sup_{a \in A}[f(x,a) + a.z]$, for some function $f$ on $\R^d\times A$, with $A$ compact set of $\R^d$, then the semi-linear PDE  \reff{semiPDE} is the Hamilton-Jacobi-Bellman equation for a stochastic control problem, where the controller can affect only the drift of the Brownian motion: $W_t + \int_0^t \alpha_s ds$,  by a  progressively measurable process $\alpha$ valued in $A$, and with a running gain function $f$.  The extension of a standard BSDE driven by a Brownian motion and an independent Poisson random  measure was considered in \cite{tanli94} and \cite{barbucpar97}, and is shown to be related in a Markovian framework to semi-linear integro-PDE.

The notion of reflected BSDE was introduced by El Karoui et al. \cite{elk97}, and consists in the addition (resp.  subtraction) of a nondecreasing process
to the standard BSDE \reff{BSDEstandard} in order to keep the solution $Y$ above (resp. below) a lower (resp. upper) obstacle, and chosen in a minimal way via the so-called Skorohod condition. Existence and uniqueness results for reflected BSDEs under  general assumptions on the obstacle have been investigated in several papers, among others \cite{ham02}, \cite{lepxu05}, \cite{pengxu05}.  We also mention works by \cite{hamouk11} and \cite{essaky} for reflected BSDEs driven by
Brownian motion and  Poisson random  measure.
An important application of reflected BSDE is its connection to optimal stopping problems and its associated variational inequalities in the Markovian case.
%which are  written  for an upper barrier in the form $U_t$ $=$ $u(t,W_t)$ as:
%\beqs
%\max \big[ - \Dt{v} - \frac{1}{2}{\rm tr}(D_x^2 v) - f^0(x,v,D_x v) ; v - u \big] &=& 0.
%\enqs

The extension  to fully nonlinear PDE, motivated in particular by uncertain volatility model and more generally by  stochastic control problem where control can affect both drift and diffusion terms of the state process, generated   important recent developments.  Soner, Touzi and Zhang \cite{sontouzha12} introduced the notion of second order BSDEs (2BSDEs), whose basic idea is to require that the solution verifies the equation
$\P^\alpha$ a.s. for every probability measure in a non dominated class of mutually
singular measures. This theory is closely related to the notion of nonlinear and $G$-expectation of Peng \cite{pen06}.
Alternatively, Kharroubi and Pham \cite{khapha12}, following \cite{khaetal10},  introduced the notion of BSDE with nonpositive jumps. The basic  idea was to constrain  the
jumps-component solution to the BSDE driven by  Brownian motion and  Poisson random  measure,  to remain nonpositive, by adding a nondecreasing process in a minimal way. A key  feature of this class of BSDEs is its formulation under a single probability measure in contrast with 2BSDEs, thus avoiding technical issues in quasi-sure analysis, and its connection with fully nonlinear HJB equation when considering a Markovian framework with a simulatable regime switching diffusion process, defined as a randomization of the controlled state process.
This approach opens new perspectives for probabilistic scheme for fully nonlinear PDEs as currently investigated in \cite{khalanpha13}.

In this paper, we define  a class of reflected BSDEs with nonpositive jumps and upper obstacle.  As in the case of doubly reflected BSDEs  with lower and upper obstacles, related to Dynkin games, our BSDE formulation involves the introduction of two nondecreasing processes, one corresponding to the nonpositive jump constraint and added in a minimal way, and the other associated to the upper reflection, satisfying the Skorohod condition, and acting in the opposite direction. The first aim of this paper is to prove the existence and uniqueness of a minimal solution to reflected BSDEs with nonpositive jumps and upper obstacle. We use a double penalization approach, and the main issue  is to obtain  uniform estimates on both penalized nondecreasing processes associated on one hand to the nonpositive jumps constraint and on the other hand to the upper obstacle. This is achieved under some regularity assumptions on the upper obstacle.  It is worth mentioning that the
running order of the  limits in the double penalization is  crucial, in contrast with the case of upper and lower reflection.  Indeed, we do not have comparison results on the jump-component solution of a BSDE, and so a priori rather few information on the sequence of nondecreasing processes associated to the jump constraint, whereas one can exploit comparison results on the $Y$-component of a BSDE in order to derive useful monotonicity property  for the sequence of nondecreasing processes associated to the upper obstacle.   Once, we get uniform estimates, we conclude by a monotonic convergence theorem for BSDEs. We also prove a dual game representation formula for the minimal solution to our BSDE, in terms of equivalent  probability measures and discount processes.

The main motivation for considering such class of  upper-reflected BSDEs with nonpositive  jumps  arises from a zero-sum stochastic differential game between a  controller and a stopper: the controller can manipulate a state process $X^\alpha$ in $\R^d$ through the selection of the control $\alpha$ valued in $A$,
while the stopper has the right to choose the duration of the game via a stopping time $\tau$. The stopper would like to minimize his expected cost:
\beq \label{coutstopper}
\E \Big[ \int_0^\tau f(X_t^\alpha,\alpha_t) dt + g(X_\tau^\alpha) \Big],
\enq
over all choices of $\tau$, while  the controller plays against him by maximizing \reff{coutstopper} over all choices of $\alpha$.  Controller-and-stopper game problem was studied in
\cite{karsudd01} when the state process $X^\alpha$ is a one-dimensional diffusion,  in \cite{karzam08} by a martingale approach and in \cite{hamlep00} by BSDE methods,  but only when the drift is controlled. General existence results for optimal actions and saddle point were recently obtained  in \cite{nutzha12} in a non Markovian and non dominated framework by exploiting the theory of
nonlinear expectations.  We also mention the recent papers \cite{matetal12}, \cite{matetal13} where the authors considered  2BSDE with reflection, in connection with
optimal stopping and Dynkin  game under nonlinear expectation.  In the Markovian case where both drift $b(X^\alpha,\alpha)$  and diffusion term $\sigma(X^\alpha,\alpha)$ of  the state process
$X^\alpha$ are controlled (hence in a non dominated framework),  the recent paper \cite{bayhua13} proved the existence of the game value, by a comparison principle for the associated Hamilton-Jacobi-Bellman Isaacs equation:
\beq
\max \Big[ - \Dt{v}  - \sup_{a\in A}\big( b(x,a).D_x v + \frac{1}{2} {\rm tr}(\sigma\sigma\trans(x,a)D_x^2 v)  +  f(x,a) \big);  & & \label{HJBintro}  \\
 v - g \big] &=& 0,  \; \mbox{ on } [0,T)\times\R^d.   \nonumber
\enq
Our second main result is to connect the minimal solution to our reflected BSDE with nonpositive jumps to a general Markovian controller-and-stopper game problem through the HJB Isaacs equation \reff{HJBintro}. We follow the idea in  \cite{bou09} and \cite{khapha12} by a randomization of the state process $X^\alpha$, and thus consider a regime switching forward diffusion process
$X$  with drift $b(X_t,I_t)$ and diffusion coefficient $\sigma(X_t,I_t)$, where $I_t$ is a pure jump process associated to the Poisson random measure driving the BSDE.  The
minimal solution $Y_t$ to the reflected BSDE with nonpositive jumps, with terminal data $\xi$ $=$ $g(X_T)$, upper obstacle $U_t$ $=$ $u(t,X_t)$, and generator $f(X_t,I_t,Y_t,Z_t)$,
is written in this Markovian  framework as: $Y_t$ $=$ $v(t,X_t,I_t)$ for some deterministic function $v$. It appears as in \cite{khapha12} that actually $v$ does not depend on $a$ in the interior of $A$ as a consequence of the non positivity jumps constraint, and we show that  $v$ is a viscosity solution to the general HJB Isaacs equation \reff{HJBintro} where the generator $f(x,a,v,\sigma\trans D_xv)$ may
depend also  on $v$ and $D_x v$.

The  rest of the paper is organized as follows. Section 2 gives a detailed formulation of reflected  BSDE with nonpositive jumps and upper obstacle. Section 3 is devoted to the existence of a minimal solution to our BSDE  by a double penalization approach.  We derive in Section 4 a dual game representation formula for the BSDE minimal solution.  Section 5 makes the connection of the minimal BSDE-solution to  fully nonlinear variational inequalities of HJB Isaacs type. We conclude in Section 6 by indicating some possible extensions to our paper.  Finally, in the appendix, we  recall some useful comparison results for BSDE with jumps, and state  a monotonic convergence theorem, which extends to the jump case the result in \cite{pengxu05}.

\section{Reflected BSDE with nonpositive jumps}

\setcounter{equation}{0} \setcounter{Assumption}{0}
\setcounter{Theorem}{0} \setcounter{Proposition}{0}
\setcounter{Corollary}{0} \setcounter{Lemma}{0}
\setcounter{Definition}{0} \setcounter{Remark}{0}

Let $(\Omega,\Fc,\P)$ be a complete probability space on which are defined  a $d$-dimensional
Brownian motion $W$ $=$ $(W_t)_{t\geq 0}$  and a Poisson random measure $\mu$ on $\R_+\times A$, where $A$ is a compact subset of $\R^q$, endowed with its Borel $\sigma$-field $\Bc(A)$.  We assume that $W$ and $\mu$ are independent, and $\mu$ has an intensity measure $\lambda(da)dt$ for some finite measure
$\lambda$  on $(A,\Bc(A))$.  We set $\tilde\mu(dt,da)$ $=$ $\mu(dt,da)-\lambda(da)dt$ the compensated martingale measure
associated to $\mu$, and denote by $\F$ $=$ $(\Fc_t)_{t\geq 0}$ the completion of the natural filtration generated by $W$ and $\mu$.

We fix a finite time duration $T<\infty$  and we denote by $\Pc$ the $\sigma$-field of $\F$-predictable subsets of $\Omega\times [0,T]$. Let us
introduce some additional notations.
We denote by:
 \begin{itemize}
 \item ${\bf L^p(}\Fc_t{\bf)}$, $p$ $\geq$ $1$, $0\leq t\leq T$, the set of $\Fc_t$-measurable random variables $X$ such that $\E|X|^p$ $<$ $\infty$.
 \item ${\bf S^2}$ the set  of real-valued c\`adl\`ag adapted processes $Y$ $=$
$(Y_t)_{0\leq t\leq T}$ such that
\[
\|Y\|_{_{{\bf S^2}}}^2 := \E\Big[ \sup_{0\leq t\leq T} |Y_t|^2 \Big] < \infty.
\]
\item ${\bf L^p(0,T)}$, $p$ $\geq$ $1$, the set of real-valued
adapted  processes  $(\phi_t)_{0\leq t\leq T}$ such that
\[
\|\phi\|_{_{\bf L^p(0,T)}}^p := \E\bigg[\int_0^T|\phi_t|^p dt\bigg] < \infty.
\]
\item ${\bf L^p(W)}$, $p$ $\geq$ $1$, the set of  $\R^d$-valued $\Pc$-measurable processes
$Z=(Z_t)_{0\leq t\leq T}$ such that
\[
\|Z\|_{_{\bf L^p(W)}}^p := \E\bigg[\bigg(\int_0^T |Z_t|^2 dt\bigg)^{\frac{p}{2}}\bigg] <\infty.
\]
\item ${\bf L^p(\tilde\mu)}$, $p$ $\geq$ $1$,  the set of
$\Pc\otimes\Bc(A)$-measurable maps $L\colon\Omega\times [0,T]\times A\rightarrow \R$ such that
\[
\|L\|_{_{{\bf L^p(\tilde\mu)}}}^p := \E\bigg[ \bigg(\int_0^T\int_A  |L_t(a)|^2 \lambda(da)dt\bigg)^{\frac{p}{2}}\bigg] < \infty.
\]
\item $\mathbf{L^2(\lambda)}$ the set of  $\Bc(A)$-measurable maps $\ell\colon A\rightarrow\R$ such that
\[
|\ell|_{_{{\bf L^2(\lambda)}}}^2 := \int_A|\ell(a)|^2\lambda(da) < \infty.
\]
\item ${\bf K^2}$  the  set of  nondecreasing predictable processes $K$ $=$ $(K_t)_{0\leq t\leq T}$ $\in$  ${\bf S^2}$ with $K_0$ $=$ $0$, so that
\beqs
\|K\|_{_{{\bf S^2}}}^2 &=&  \E|K_T|^2.
\enqs
\end{itemize}

We are then  given three  objects:
\begin{itemize}
\item[1.] A {\it terminal  condition} $\xi$ $\in$ ${\bf L^2(}\Fc_T{\bf)}$.
%, i.e., a real-valued  $\Fc_T$-measurable random variable, which is square-integrable:
%\beqs \label{square2xi}
%\E\big[|\xi|^2\big] &<& \infty.
%\enqs
\item[2.] A  {\it generator function}
$F:\Omega\times[0,T]\times\R\times\R^d \times \mathbf{L^2(\lambda)} \rightarrow \R$,
which is a $\Pc\otimes\Bc(\R)\otimes\Bc(\R^d)\otimes\Bc(\mathbf{L^2(\lambda)})$-measurable map, satisfying:
\begin{itemize}
\item[(i)] The square integrability condition:
\beqs
\label{IntCondF}
\E\bigg[ \int_{0}^T |F(t,0,0,0)|^2 dt\bigg] & < & \infty.
\enqs
\item[(ii)] The uniform Lipschitz condition:
%there exists a constant $C_F$ such that
\beqs \label{lipF}
|F(t,y,z,\ell)-F(t,y',z',\ell')| & \leq & C_F \big(|y-y'|+|z-z'| + |\ell-\ell'|_{_{{\bf L^2(\lambda)}}} \big)\;,
\enqs
for all $t \in [0,T]$, $y,y' \in \R$, $z,z'\in \R^d$, and $\ell,\ell'\in \mathbf{L^2(\lambda)}$, where $C_F$ is some positive constant.
\item[(iii)] The monotonicity condition:
\beq
\hspace{-1cm}F(t,y,z,\ell)-F(t,y,z,\ell') & \leq & \int_A (\ell(a)-\ell'(a))\gamma(t,y,z,\ell,\ell',a) \lambda(da),  \label{monol}
\enq
for all  $t \in[0,T]$, $y \in \R$, $z\in \R^d$,  and $\ell,\ell'\in \mathbf{L^2(\lambda)}$, where  $\gamma$ $:$  $\Omega\times [0,T]\times\R\times\R^d\times\mathbf{L^2(\lambda)}\times\mathbf{L^2(\lambda)}\times A \rightarrow \R$  is a
$\Pc\otimes\Bc(\R)\otimes\Bc(\R^d)\otimes\Bc(\mathbf{L^2(\lambda)})\otimes\Bc(\mathbf{L^2(\lambda)})\otimes\Bc(A)$-measurable map satisfying: $0 \leq \gamma(t,y,z,\ell,\ell',a) \leq C_\gamma$,
for all $t \in [0,T]$, $y \in \R$, $z \in \R^d$, $\ell,\ell'\in \mathbf{L^2(\lambda)}$, and $a$ $\in$ $A$, for some positive constant $C_\gamma$.
\end{itemize}
%\item[3.]  A {\it lower jump constraint function} $h$ $:$ $\Omega\times [0,T]\times\R\times A$ $\rightarrow$ $\R$, given by:
%\beqs
%h(t,\ell,a) = \ell - c_t(a),
%\enqs
%where $c$ $:$ $\Omega\times[0,T]\times A$ $\rightarrow$ $\R$ is a
%$\Pc\otimes\Bc(A)$-measurable map, such that
%\beqs
%\E\bigg[ \int_0^T \int_A |c_t(a)|^2 \lambda(da) dt \bigg] &<& \infty
%\enqs
%and $c_t(a)$ $\geq$ $0$, $d\P\otimes dt\otimes\lambda(da)$ a.e..
\item[3.] An {\it upper barrier} $U$ $\in$ ${\bf S^2}$ satisfying $U_T$ $\geq$ $\xi$, almost surely.
\end{itemize}

\vspace{2mm}

Let us now consider our problem of  reflected BSDE with nonpositive jumps. We say that a quintuple $(Y,Z,L,K^+,K^-)$ $\in$
${\bf S^2}\times{\bf L^2(W)}\times{\bf L^2(\tilde \mu)}\times{\bf K^2}\times{\bf K^2}$ is a solution to the upper-reflected BSDE with nonpositive jumps with data $(\xi,F,U)$ if the following relation holds:
\beq
Y_t & = & \xi + \int_t^T F(s,Y_s,Z_s,L_s) ds    +  K_T^+ - K_t^+ - (K_T^--K_t^-)   \label{BSDEgen}  \\
& &  -  \int_t^T  Z_s dW_s  -  \int_t^T\int_A   L_s(a)  \mu(ds,da), \;\;\;\;\;   0 \leq t \leq T,  \; a.s. \nonumber
\enq
together with the jump constraint
\beq
L_t(a) & \leq & 0\;, \;\;\;\;\; d\P\otimes dt\otimes\lambda(da) \;\;  a.e. \label{Lcons}
\enq
and the upper constraint
\beq
Y_t & \leq & U_t\;, \;\;\;\;\; 0 \leq t \leq T,\;a.s. \label{Ucons} \\
\int_0^T (U_{t^-} - Y_{t^-}) dK^-_t & = & 0\;, \;\;\;\;\;\;\, a.s. \label{Sko}
\enq

We look for the \emph{minimal solution} $(Y,Z,L,K^+,K^-)$, in the sense that for any other solution $(\tilde Y, \tilde Z,\tilde L,\tilde K^+,\tilde K^-)$ to the reflected BSDE with nonpositive jumps \reff{BSDEgen}-\reff{Lcons}-\reff{Ucons}-\reff{Sko}, it must hold that $Y$ $\leq$ $\tilde Y$.

\begin{Remark}
{\rm
We have chosen to formulate the BSDE \reff{BSDEgen}  directly in terms of the random measure $\mu$ instead of the compensated random measure $\tilde\mu$ since we dealt with
finite  intensity measure $\lambda(A)$ $<$ $\infty$. Of course, one can formulate equivalently the BSDE \reff{BSDEgen} in terms of $\tilde\mu$ by changing  the generator $F$ to:
\beqs
\tilde F(t,y,z,\ell) &=& F(t,y,z,\ell) - \int_A \ell(a) \lambda(da).
\enqs
In this case, the monotonicity condition \reff{monol} for $\tilde F$ holds with a measurable map $\tilde\gamma$ satisfying:
$-1 \leq \tilde\gamma(t,y,z,\ell,\ell',a) \leq C_{\tilde\gamma}$,  for all $t \in [0,T]$, $y \in \R$, $z \in \R^d$, $\ell,\ell'\in \mathbf{L^2(\lambda)}$, and $a$ $\in$ $A$, for some positive constant $C_{\tilde\gamma}$.
\ep
}
\end{Remark}

\begin{Remark}
\label{RemarkUniq}
Uniqueness of the minimal solution. {\rm Uniqueness of a  minimal solution holds in the following sense: if
$(Y,Z,L,K^+,K^-)$ and $(Y,\tilde Z,\tilde L,\tilde K^+,\tilde K^-)$ are minimal solutions to \reff{BSDEgen}-\reff{Lcons}-\reff{Ucons}-\reff{Sko}, then
$Y$ $=$ $Y'$, $Z$ $=$ $Z'$, $L$ $=$ $L'$, and $K^+-K^-$ $=$ $\tilde K^+-\tilde K^-$.   As a matter of fact, the uniqueness of the $Y$ component is clear by definition. Then, denoting by $K:$ $=$ $K^+-K^-$, and $\tilde K$ $:=$ $\tilde K^+-\tilde K^-$, which are predictable finite variation processes,
we have
\beqs
\int_0^t \big[F(s,Y_s,Z_s,L_s) - F(s,Y_s,\tilde Z_s,\tilde L_s)\big] ds + K_t -  \tilde K_t & & \\
\; + \; \int_0^t (\tilde Z_s-Z_s)dW_s +
\int_0^t \int_A (\tilde L_s(a)-L_s(a))\mu(ds,da) &=& 0,
\enqs
for all $t$ $\in$ $[0,T]$, almost surely. The uniqueness of $Z$ $=$ $\tilde Z$ follows by identifying the Brownian part and the finite variation part, while the uniqueness of $(L,K)$ $=$ $(\tilde L,\tilde K)$  is obtained by identifying the predictable part, and by recalling that the jumps of $\mu$ are totally inaccessible.
\ep
%we see by taking the predictable projection in the above relation, and identifying Brownian parts and finite variation parts that: $L$ $=$ $\tilde L$ and $Z$ $=$ $\tilde Z$.  Consequently, we get $K$ $=$ $\tilde K$.
%Actually, uniqueness of $K^+$ and $K^-$ hold when they are the Jordan components of the  and by uniqueness of the Jordan decomposition (see, e.g., %Lemma 1.35 in \cite{jac79}), we conclude that $K^+$ $=$ $\tilde K^+$ and $K^-$ $=$ $\tilde K^-$.
}
\end{Remark}

\vspace{2mm}

The main feature in this class of BSDEs is to consider a reflection constraint on $Y$ in addition to the nonpositive jump constraint as already studied in \cite{khaetal10} and \cite{khapha12}. Moreover,  we deal with  an upper barrier $U$  associated to a nondecreasing process $K^-$, which is subtracted in \reff{BSDEgen}  from the nondecreasing process $K^+$ associated to the nonpositive constrained jumps.
In order to ensure that the problem of getting a minimal solution to \reff{BSDEgen}-\reff{Lcons}-\reff{Ucons}-\reff{Sko}  is well-posed, and  similarly as in \cite{khapha12}, we make the assumption that there exists a
supersolution to the BSDE with nonpositive jumps, namely:

\vspace{3mm}

\noindent {\bf (H0)}  \hspace{7mm}  There exists $(\bar Y,\bar Z,\bar L, \bar K^+)$ $\in$ ${\bf S^2}\times{\bf L^2(W)}\times{\bf L^2(\tilde \mu)}\times{\bf K^2}$ satisfying the BSDE with nonpositive jumps:
\beq
\bar Y_t & = & \xi + \int_t^T F(s,\bar Y_s, \bar Z_s, \bar L_s) ds    + \bar K_T^+ - \bar K_t^+ \label{H01} \\
& &  -  \int_t^T  \bar Z_s dW_s  -  \int_t^T\int_A  \bar L_s(a)  \mu(ds,da)\;, \;\;\;\;\;   0 \leq t \leq T,  \; a.s. \nonumber
\enq
and
\beq
\bar L_t(a)  & \leq & 0\;, \;\;\;\;\; d\P\otimes dt\otimes\lambda(da) \;\;  a.e. \label{H02}
\enq
We shall see later in the Markovian case (see Remark \ref{remH0})  how this condition {\bf (H0)} is directly satisfied.

\section{Existence and approximation by double penalization}

\setcounter{equation}{0} \setcounter{Assumption}{0}
\setcounter{Theorem}{0} \setcounter{Proposition}{0}
\setcounter{Corollary}{0} \setcounter{Lemma}{0}
\setcounter{Definition}{0} \setcounter{Remark}{0}

This section is devoted to the existence of the minimal solution to \reff{BSDEgen}-\reff{Lcons}-\reff{Ucons}-\reff{Sko}.   We use a penalization approach and introduce the doubly indexed sequence of BSDEs with jumps:
\beq
Y_t^{n,m} &=&   \xi + \int_t^T F(s,Y_s^{n,m},Z_s^{n,m},L_s^{n,m}) ds + K_T^{n,m,+} - K_t^{n,m,+} - (K_T^{n,m,-}-K_t^{n,m,-}) \nonumber \\
& & -  \int_t^T  Z_s^{n,m} dW_s  -  \int_t^T\int_A  L_s^{n,m}(a)  \mu(ds,da),  \label{penBSDE}
\enq
for $n,m$ $\in$ $\N$, where $K^{n,m,+}$ and $K^{n,m,-}$ are the nondecreasing continuous processes in ${\bf K^2}$ defined by
\[
K_t^{n,m,+} \; = \;  m \int_0^t \int_A (L_s^{n,m}(a))_+ \lambda(da) ds, \qquad
K_t^{n,m,-} \; = \;  n \int_0^t (U_s - Y_s^{n,m})_- ds.
\]
Here we use the notation $f_+$ $=$ $\max(f,0)$ and $f_-$ $=$ $\max(-f,0)$ to denote the positive and negative parts of $f$.  Notice that this penalized BSDE can be written as
\[
Y_t^{n,m} = \xi + \int_t^T F_{n,m}(s,Y_s^{n,m},Z_s^{n,m},L_s^{n,m}) ds - \int_t^T  Z_s^{n,m} dW_s - \int_t^T\int_A   L_s^{n,m}(a)  \mu(ds,da),
\]
with a generator $F_{n,m}$ given by
\beqs
F_{n,m}(t,y,z,\ell) &=& F(t,y,z,\ell) + m \int_A (\ell(a))_+ \lambda(da) -  n (U_t - y)_-, \;\;\;  a.s.
\enqs
for $(t,y,z,\ell)$ $\in$ $[0,T]\times\R\times\R^d\times\mathbf{L^2(\lambda)}$. Observe that the generator $F_{n,m}$ satisfies the assumptions of square integrability and  uniform Lipschitzianity, which ensure by Lemma 2.4 in \cite{tanli94} the existence and uniqueness of a solution $(Y^{n,m},Z^{n,m},L^{n,m})$ $\in$ ${\bf S^2}\times{\bf L^2(W)}\times{\bf L^2(\tilde \mu)}$  to the  BSDE with jumps \reff{penBSDE}. Notice also that $F_{n,m}$ satisfies the monotonicity condition
\reff{monol},  is increasing in $m$ for any  fixed $n$,  and decreasing in $n$ for any fixed $m$. Thus,  by the comparison Theorem \ref{CompThm}, we deduce that $(Y^{n,m})_{n,m}$ inherits the same property:
\beq \label{Ynmcroidecroi}
Y^{n+1,m} & \leq & Y^{n,m} \;\, \leq \;\, Y^{n,m+1}, \;\;\;\; \forall\, n,m \in \N.
\enq
We shall first fix $m$, and let $n$ to infinity, and then let $m$ to infinity (the order of the limits is important here, see
Remark \ref{limitorder}).
The key point, as in the case of doubly reflected BSDEs related to Dynkin games, is
to deal with the difference of  the nondecreasing processes $K^{n,m+}$ and $K^{n,m,-}$, and the main difficulty is to prove their  convergence  towards respectively the  nondecreasing
processes $K^+$ and $K^-$, which appear in the minimal solution to the reflected BSDE with nonpositive jumps we are looking for. We have  to impose some regularity conditions on the upper barrier process that will be precised later.

\vspace{3mm}

For fixed $m$, let us now consider the reflected BSDE with jumps:
\beq
Y_t^m & = & \xi + \int_t^T F_m(s,Y_s^m,Z_s^m,L_s^m) ds    - (K_T^{m,-} - K_t^{m,-})   \label{BSDEmaxm} \\
 & &  -  \int_t^T Z_s^m dW_s  -  \int_t^T\int_A L_s^m(a) \mu(ds,da)\;, \;\;\;\;\;   0 \leq t \leq T,  \; a.s. \nonumber
\enq
and
\beq
Y_t^m & \leq & U_t\;, \;\;\;\;\; 0 \leq t \leq T,\;a.s. \label{consUm} \\
\int_0^T (U_{t^-} - Y_{t^-}^m) dK_t^{m,-} & = & 0\;, \;\;\;\;\;\;\, a.s. \label{Skom}
\enq
where
\beq \label{defFm}
F_m(t,y,z,\ell) &=& F(t,y,z,\ell) + m \int_A (\ell(a))_+ \lambda(da)\;, \;\;\;\;\; a.s.
\enq
for $(t,y,z,\ell)$ $\in$ $[0,T]\times\R\times\R^d\times\mathbf{L^2(\lambda)}$. We know from Theorem 4.2 in \cite{hamouk11} that there exists a unique solution
$(Y^m,Z^m,L^m,K^{m,-})$ $\in$ ${\bf S^2}\times{\bf L^2(W)}\times{\bf L^2(\tilde \mu)}\times{\bf K^2}$ to the reflected BSDE with jumps \reff{BSDEmaxm}-\reff{consUm}-\reff{Skom}.

\begin{Remark}
\label{R:HamEssaky}
{\rm Note that in \cite{hamouk11} the existence of $(Y^m,Z^m,L^m,K^{m,-})$ is proved using a fixed point argument and not through the penalized sequence $(Y^{n,m},Z^{n,m},L^{n,m})$, except for the particular case where the generator $F_{n,m}(t,\omega)$ does not depend on $y,z,\ell$, see Theorem 4.1 and Remark 4.1(i) in \cite{hamouk11}. The reason is that in \cite{hamouk11} the authors do not impose any monotonicity condition on the generator $F$ and therefore they do not have at disposal a comparison theorem for BSDEs with jumps. Nevertheless, under our monotonicity condition \eqref{monol} and by means of the comparison Theorem \ref{CompThm}, the existence of $(Y^m,Z^m,L^m,K^{m,-})$ can be proved via the penalized sequence $(Y^{n,m},Z^{n,m},L^{n,m})$. This program is carried out in \cite{essaky}, Theorem 5.1, even though under the additional hypothesis that the barrier $U$ is a $\Pc$-measurable process. More precisely, it can be shown that $Y^m$ is obtained as the decreasing limit of $Y^{n,m}$ when $n$ goes to infinity:
\beqs
Y_t^m &=& \lim_{n\rightarrow\infty}\downarrow Y_t^{n,m}, \;\;\; 0 \leq t \leq T, a.s.
\enqs
and this convergence also holds in ${\bf L^2(0,T)}$. Furthermore, $(Z^{n,m},L^{n,m})$ converges weakly to  $(Z^m,L^m)$  in  ${{\bf L^2(W)}\times \bf L^2(\tilde \mu)}$,  and we have the strong convergence
\beqs
(Z^{n,m},L^{n,m})  & \rightarrow & (Z^m,L^m)  \;\; \mbox{ in }  \;  {{\bf L^p(W)}\times \bf L^p(\tilde \mu)}, \;\;\; \mbox{ as } \; n \rightarrow \infty,
\enqs
for any $p$ $\in$ $[1,2)$, while
\beqs
K_t^{n,m,-}  & \rightharpoonup & K_t^{m,-}  \;\;\; \mbox{ weakly in }  {\bf L^2}(\Fc_t),  \;\;\; \mbox{ as } \; n \rightarrow \infty
\enqs
for all $0 \leq t \leq T$.
\ep
}
\end{Remark}

\vspace{2mm}

We first derive the following important  property on the sequence of nondecreasing processes $(K^{m,-})$.

\begin{Lemma} \label{lemKm}
The sequence of processes $(K^{m,-})_m$ satisfies:
\beq
\label{K-croissant}
K_t^{m,-} - K_s^{m,-} & \leq & K_t^{m+1,-}  -  K_s^{m+1,-}, \;\;\;\;\; 0 \leq s \leq t \leq T,\;a.s., \;\; \forall m \in \N.
\enq
\end{Lemma}
{\bf Proof.} By definition of $K^{n,m,-}$, and from \reff{Ynmcroidecroi}, we clearly have for all $n,m$ $\in$ $\N$:
\beqs
K_t^{n,m,-} - K_s^{n,m,-} &\leq& K_t^{n,m+1,-} - K_s^{n,m+1,-}, \;\;\;\;\; 0 \leq s \leq t \leq T,\;a.s.
\enqs
Thus, by passing to the (weak) limit as $n$ goes to infinity, we get the required result.
\ep

\vspace{3mm}

%Let us also introduce the nondecreasing process in ${\bf K^2}$:
%\beq
%\label{Km+}
%K_t^{m,+} &=&  m \int_0^t \int_A h^+(L_s^{m}(a)) \lambda(da) ds\;, \;\;\;\;\; 0 \leq t\leq T.
%\enq
By \reff{Ynmcroidecroi}, we see that $(Y^m)_m$ is a nondecreasing sequence: $Y^{m}$ $\leq$ $Y^{m+1}$, and we denote:
\beqs
\underline Y_t &:=& Y_t^0, \;\;\;\;\; 0 \leq t \leq T,   \label{underlineY}
\enqs
which thus provides a lower bound for the sequences $(Y^m)$ and $(Y^{n,m})$:
\beq  \label{lowerYnm}
\underline Y_t \; \leq \; Y_t^m  & \leq & Y_t^{n,m}, \;\;\; 0 \leq t \leq T, \; \forall n,m \in \N.
\enq
Moreover, under condition {\bf (H0)}, we observe  that the quintuple  $(\bar Y,\bar Z,\bar L,\bar K^+,\bar K^-)$ satisfies $\int_A (\bar L_t(a))_+ \lambda(da)$ $=$ $0$ $dt\otimes d\P$ a.e. so that
\beqs
F_{n,m}(t,\bar Y_t,\bar Z_t,\bar L_t) & \leq & F(\bar Y_t,\bar Z_t,\bar L_t), \;\;\; dt\otimes d\P \; a.e.
\enqs
By the comparison Theorem~\ref{CompThm}, we  then get an upper bound for the sequences $(Y^m)$ and  $(Y^{n,m})$:
\beq \label{upperYnm}
Y_t ^m \; \leq \; Y_t^{n,m} & \leq & \bar Y_t, \;\;\;  0 \leq t \leq T, \forall n,m \in \N.
\enq

\vspace{3mm}

By standard arguments, we now state some estimates on the doubly indexed sequence $(Y^{n,m},Z^{n,m},L^{n,m},K^{n,m,+})$ expressed in terms of  $(K^{n,m,-})$.

\begin{Lemma}
\label{L:ZLbdd}
Let assumption {\bf (H0)} hold. Then there exists a positive constant $C$, such that for all  $n,m$ $\in$ $\N$,
\beq
 & & \|Y^{n,m}\| _{_{\bf S^2}}^2 +   \|Z^{n,m}\|_{_{\bf L^2(W)}}^2 + \|L^{n,m}\|_{_{\bf L^2(\tilde \mu)}}^2 + \|K^{n,m,+}\|_{_{\bf S^2}}^2  \nonumber \\
 & \leq&  C \bigg( \E|\xi|^2 + \E\int_0^T |F(s,0,0,0)|^2 ds + \big\| \underline Y \big\|_{_{\bf S^2}}^2  + \big\| \bar Y \big\|_{_{\bf S^2}}^2 + \|K^{n,m,-}\|_{_{\bf S^2}}^2\bigg).  \label{estimY}
\enq
\end{Lemma}
{\bf Proof.}
In what follows we shall denote by $C$ $>$ $0$ a generic positive constant depending only on $T$, $\lambda(A)$, and the Lipschitz constant of $F$, which may vary from line to line.
Proceeding as in the proof of Lemma 3.3 in \cite{khapha12}, we apply It\^o's formula to $|Y_s^{n,m}|^2$ between $t$ and $T$, and get after some rearrangement:
\beq
& & \E|Y_t^{n,m}|^2 + \|Z^{n,m}1_{[t,T]}\|_{_{\bf L^2(W)}}^2+ \|L^{n,m}1_{[t,T]}\|_{_{\bf L^2(\tilde \mu)}}^2  \nonumber \\
&=& \E|\xi|^2  +  2\E\int_t^T Y_s^{n,m}F(s,Y_s^{n,m},Z_s^{n,m},L_s^{n,m})ds  -  2\E\int_t^T\int_A Y_{s^-}^{n,m}L_s^{n,m}(a) \lambda(da)ds \nonumber \\
& & \; + \; 2\E\int_t^T Y_s^{n,m}dK_s^{n,m,+} -  2\E\int_t^T Y_s^{n,m}dK_s^{n,m,-}.  \label{Ynmito}
\enq
By the linear growth condition on  $F$, the inequality $ab$ $\leq$ $a^2/2$ $+$ $b^2/2$, and recalling that $\lambda(A)$ $<$ $\infty$,  we get
\beq
& & \!\!\!\!\! 2\E\int_t^T Y_s^{n,m}F(s,Y_s^{n,m},Z_s^{n,m},L_s^{n,m})ds - 2\E\int_t^T\int_A Y_{s^-}^{n,m}L_s^{n,m}(a) \lambda(da)ds   \label{YnmF} \\
& \leq & \!\!\!\!\! C\E\int_t^T|Y_s^{n,m}|^2 ds + \frac{1}{2}\E\int_0^T |F(s,0,0,0)|^2 ds + \frac{1}{2}\|Z^{n,m}1_{[t,T]}\|_{_{\bf L^2(W)}}^2  + \frac{1}{2}\|L^{n,m}1_{[t,T]}\|_{_{\bf L^2(\tilde \mu)}}^2. \nonumber
\enq
From the bounds \reff{lowerYnm}-\reff{upperYnm} on $Y^{n,m}$:  $\underline Y$ $\leq$ $Y^{n,m}$ $\leq$ $\bar Y$,
and thanks to the inequality $2ab$ $\leq$ $a^2/\alpha$ $+$ $\alpha b^2$ for any constant $\alpha>0$, we have
\beqs
& & 2\E\int_t^T Y_s^{n,m}dK_s^{n,m,+}  - 2\E\int_t^T Y_s^{n,m}dK_s^{n,m,-} \nonumber \\
&\leq&  \frac{1}{\alpha} \Big( \big\| \underline Y \big\|_{_{\bf S^2}}^2 + \big\| \bar Y \big\|_{_{\bf S^2}}^2 \Big) + \alpha \E|K_T^{n,m,+}-K_t^{n,m,+}|^2 + \alpha \E|K_T^{n,m,-}-K_t^{n,m,-}|^2 \nonumber \\
& \leq & \frac{1}{\alpha} \Big( \big\| \underline Y \big\|_{_{\bf S^2}}^2 + \big\| \bar Y \big\|_{_{\bf S^2}}^2 \Big) + 3\alpha \E|K_T^{n,m,-}-K_t^{n,m,-}|^2 + 2\alpha \E|K_T^{n,m}-K_t^{n,m}|^2,
\enqs
where we set  $K_t^{n,m}$ $:=$ $K_t^{n,m,+}$ $-$ $K_t^{n,m,-}$, so that $\E|K_T^{n,m,+}-K_t^{n,m,+}|^2$ $\leq$ $2\E|K_T^{n,m}-K_t^{n,m}|^2$ $+$ $2\E|K_T^{n,m,-}-K_t^{n,m,-}|^2$. Together with \reff{YnmF} and
\reff{Ynmito}, this yields:
\beq
 & & \E|Y_t^{n,m}|^2  +   \frac{1}{2} \|Z^{n,m}1_{[t,T]}\|_{_{\bf L^2(W)}}^2 +  \frac{1}{2} \|L^{n,m}1_{[t,T]}\|_{_{\bf L^2(\tilde \mu)}}^2  \nonumber \\
 & \leq &  C\E\int_t^T|Y_s^{n,m}|^2 ds +   \E|\xi|^2  +  \frac{1}{2}\E\int_0^T |F(s,0,0,0)|^2 ds
 + \frac{1}{\alpha} \Big( \big\| \underline Y \big\|_{_{\bf S^2}}^2 + \big\| \bar Y \big\|_{_{\bf S^2}}^2 \Big) \nonumber \\ 
 & & \;\;\;  + \;  3\alpha \E|K_T^{n,m,-}-K_t^{n,m,-}|^2 + 2\alpha \E|K_T^{n,m}-K_t^{n,m}|^2.  \label{YnmK}
\enq
Now, from the relation \reff{penBSDE}, we have
\beqs
K_T^{n,m} - K_t^{n,m} &=& Y_t^{n,m} - \xi - \int_t^T F(s,Y_s^{n,m},Z_s^{n,m},L_s^{n,m}) ds   \\
& & \; + \int_t^T Z_s^{n,m} dW_s + \int_t^T\int_A L_s^{n,m}(a) \mu(ds,da),
\enqs
so that by the linear growth condition on $F$:
\beq
\E|K_T^{n,m}-K_t^{n,m}|^2 &\leq& C \bigg(   \E|\xi|^2 + \E\int_0^T |F(s,0,0,0)|^2 ds  +  \E|Y_t^{n,m}|^2  \label{estimK} \\
&& \;\;\;\;\;  + \;  \E\int_t^T |Y_s^{n,m}|^2 ds +  \|Z^{n,m}1_{[t,T]}\|_{\bf L^2(W)}^2 + \|L^{n,m}1_{[t,T]}\|_{\bf L^2(\tilde \mu)}^2 \bigg). \nonumber
%\\
%& \leq &  C \Big(   \E|\xi|^2 + \E\int_0^T |F(s,0,0,0)|^2 ds  +  \big\| \underline Y \big\|_{_{\bf S^2}}^2 + \big\| \bar Y \big\|_{_{\bf S^2}}^2 \\
%&& \;\;\;\;\;  + \;  \E\int_t^T |Y_s^{n,m}|^2 ds +  \|Z^{n,m}\|_{\bf L^2(W)}^2 + \|L^{n,m}\|_{\bf L^2(\tilde \mu)}^2 \Big),
\enq
%where we used  again the bound: $\underline Y$ $\leq$ $Y^{n,m}$ $\leq$ $\bar Y$.
By choosing $\alpha>0$ such that $2\alpha C$ $\leq$ $1/4$, and plugging this estimate of  $\E|K_T^{n,m}-K_t^{n,m}|^2$ into
\reff{YnmK},  we get for all $0\leq t\leq T$:
\beq
& & \frac{3}{4} \E|Y_t^{n,m}|^2  +   \frac{1}{4} \|Z^{n,m}1_{[t,T]}\|_{_{\bf L^2(W)}}^2 +  \frac{1}{4} \|L^{n,m}1_{[t,T]}\|_{_{\bf L^2(\tilde \mu)}}^2 \nonumber \\
& \leq &  C\E\int_t^T|Y_s^{n,m}|^2 ds +   \frac{5}{4} \E|\xi|^2  +  \frac{3}{4}\E\int_0^T |F(s,0,0,0)|^2 ds   \nonumber \\
& & \; + \;  \frac{1}{\alpha} \Big( \big\| \underline Y \big\|_{_{\bf S^2}}^2 + \big\| \bar Y \big\|_{_{\bf S^2}}^2 \Big) + 3\alpha \E|K_T^{n,m,-}-K_t^{n,m,-}|^2 \nonumber \\
& \leq &  C\bigg( \big\| \underline Y \big\|_{_{\bf S^2}}^2 + \big\| \bar Y \big\|_{_{\bf S^2}}^2   +   \E|\xi|^2  +  \E\int_0^T |F(s,0,0,0)|^2 ds \bigg)
+ 12\alpha \|K^{n,m,-}\|_{\bf S^2}^2,  \label{interL}
\enq
where we used again the bounds $\underline Y\leq Y^{n,m}\leq \bar Y$ and the inequality $\E|K_T^{n,m,-}-K_t^{n,m,-}|^2$ $\leq$ $4\E|K_T^{n,m,-}|^2$. This proves, taking $t=0$ in \eqref{interL}, the required  estimate \reff{estimY} for $(Z^{n,m},L^{n,m})$, and also for $K^{n,m,+}$ by \reff{estimK}, and recalling that $\E|K_T^{n,m,+}|^2$ $\leq$ $2\E|K_T^{n,m}|^2$ $+$ $2\E|K_T^{n,m,-}|^2$. Finally, the estimate for $\|Y^{n,m}\| _{_{\bf S^2}}$ in \reff{estimY}  follows as usual from the relation \reff{penBSDE}, Burkholder-Davis-Gundy inequality, and the  estimates for $(Z^{n,m},L^{n,m},K^{n,m,+})$.
\ep

\vspace{3mm}

The key point is now to obtain a uniform estimate on $K^{n,m,-}$,  and consequently  uniform estimates on $(Y^{n,m},Z^{n,m},L^{n,m},K^{n,m,+})$ in view of  Lemma \ref{L:ZLbdd}.
Let us introduce the  following set of probability measures.  For $m$ $\in$ $\N$, let $\Vc_m$ be the set of $\Pc\otimes\Bc(A)$-measurable processes valued in $(0,m]$,  $\Vc$ $=$ $\cup_m \Vc_m$, and given $\nu$ $\in$ $\Vc$, consider the probability measure $\P^\nu$ equivalent to $\P$ on $(\Omega,\Fc_T)$ with Radon-Nikodym density:
\beqs
\frac{d\P^\nu}{d\P} \bigg|_{\Fc_t}  &=& \zeta_t^\nu \; := \;  \Ec_t \bigg(\int_0^. \int_A (\nu_s(a) - 1) \tilde\mu(ds,da) \bigg),
\enqs
where $\Ec_t(\cdot)$ is the Dol\'eans-Dade exponential. Indeed, since $\nu$ $\in$ $\Vc$ is essentially bounded, and $\lambda(A)$ $<$ $\infty$, it is known that  $\zeta^\nu$ is a uniformly integrable martingale (see e.g. Lemma 4.1 in \cite{khapha12}), and so defines  a probability measure $\P^\nu$. Moreover, $\zeta_T^\nu$ $\in$ ${\bf L^p}(\Fc_T)$ for any $p$ $\geq$ $1$.
Notice that the Brownian motion $W$ remains a Brownian motion $W$ under $\P^\nu$, while the effect of the probability measure $\P^\nu$, by Girsanov's theorem, is to change the
compensator $\lambda(da)dt$ of $\mu$ under $\P$ to $\nu_t(a)\lambda(da)dt$ under $\P^\nu$.  We then denote by $\tilde\mu^\nu(dt,da)$ $:=$ $\mu(dt,da)-\nu_t(a)\lambda(da)dt$ the compensated martingale measure of $\mu$ under $\P^\nu$.

Inspired by \cite{hamlepmat97} (see also \cite{cvikar96}), we make the following regularity  assumption on the upper barrier:

\vspace{3mm}

\noindent {\bf (H1)}  \hspace{3mm}  There exists a nonincreasing sequence of processes $(U^k)_k$ such that:
\begin{itemize}
\item[(i)] $\lim_{k\rightarrow\infty}U_t^k$ $=$ $U_t$, for all $0$ $\leq$ $t$ $\leq$ $T$, a.s..
\item[(ii)] For any $k\in\N$, $U^k$ is in the form:
\beqs
U_t^k &=& U_0^k + \int_0^t \upsilon_s^k ds + \int_0^t \vartheta_s^k dW_s, \;\;\;\;\; 0 \leq t \leq T, \; a.s.
\enqs
where $(\ups^k)_k\subset{\bf L^2(0,T)}$ and $(\var^k)_k\subset{\bf L^2(W)}$.
\item[(iii)] There exists  some $p$ $>$ $2$  such that:
\beqs
\sup_{k\in\N} \int_0^T \E\Big[  \esssup_{\nu\in\Vc} \E^\nu\big[\sup_{t \leq s \leq T}\big(|U_s^k|^p + |\upsilon_s^k|^p + |\vartheta_s^k|^p\big) \big| \Fc_t \big] \Big]dt  & &  \nonumber \\
\; + \; \int_0^T\E\Big[\esssup_{\nu\in\Vc}\E^\nu\big[\sup_{t \leq s \leq T}\big|F(s,0,0,0)\big|^p \big| \Fc_t \big]\Big] dt &  < & \infty.  \label{IntCondU}
\enqs
\end{itemize}

We shall see later  in the Markovian framework how Assumption {\bf (H1)} is automatically satisfied, see Remark \ref{R:H1Markovian}.  The following key lemma  states a uniform estimate for $K^{n,m,-}$ under condition {\bf (H1)}.

 \vspace{3mm}

\begin{Lemma} \label{lemKnm-}
Under condition {\bf (H1)}, we have
\beqs
\sup_{n,m\in\N}  \big\| K^{n,m,-} \big\|_{_{\bf S^2}} &< & \infty.
\enqs
\end{Lemma}
{\bf Proof.}
%\marginpar{I don't  think we need this condition on $F$ with the monotonicity condition}
%\beq
%\label{AssumptionF1}
%|F(t,0,0,\ell)| &\leq& \bar C\big(1+|F(t,0,0,0)|\big), \qquad (t,\ell)\in[0,T]\times\mathbf{L^2(\lambda)},\;a.s.
%\enq
%and
%\beq
%\label{AssumptionF2}
%\int_0^T\E\Big[\esssup_{\nu\in\Vc_m}\E^\nu\Big[\sup_{t \leq s \leq T}\big|F(s,0,0,0)\big|^p \Big| \Fc_t \Big]\Big] dt &\leq& \bar C
%\enq
%for some $p>2$ such that \reff{IntCondU} holds.
Let $(U^k)_k$ be in the form  as in assumption {\bf (H1)}(ii) and consider for  positive integers $n,m,k$,
the difference $\bar Y^{n,m,k}$ $:=$ $Y^{n,m} - U^k$,    which is then expressed in backward form as:
\beq
\bar Y_t^{n,m,k}  &=&   \xi - U_T^k + \int_t^T \big(F(s,Y_s^{n,m},Z_s^{n,m},L_s^{n,m}) +  \ups_s^k\big) ds  \nonumber \\
& &  \; +\; m\int_t^T\int_A (L_s^{n,m}(a))_+\lambda(da)ds   - n\int_t^T(U_s - U_s^k -  \bar Y_s^{n,m,k})_-ds  \nonumber \\
& &  \; - \;  \int_t^T  \big(Z_s^{n,m} -  \var_s^k\big) dW_s  -  \int_t^T\int_A  L_s^{n,m}(a) \mu(ds,da). \label{barYnm}
\enq
Now, by the Lipschitz condition of $F$ in $(y,z)$, and the monotonicity condition \reff{monol} of $F$ in $\ell$, we have for all $n,m$ $\in$ $\mathbb{N}$:
\beqs
F(t,Y_t^{n,m},Z_t^{n,m},L_t^{n,m}) &=& F(t,0,0,0) +  \alpha_t^{n,m} Y_t^{n,m} + \beta_t^{n,m}.  Z_t^{n,m}  \\
& & \;\;\;  + \;  \int_A \gamma_t^{n,m}(a) L_t^{n,m}(a)  \lambda(da) -  \delta_t^{n,m},
\enqs
for some  sequence of bounded predictable processes $(\alpha^{n,m})$ valued in $\R$,  $(\beta^{n,m})$ valued in $\R^d$,  uniformly bounded in $n,m$,
a nonnegative sequence of predictable process $(\delta^{n,m})$, and a nonnegative sequence of  bounded $\Pc\otimes\Bc(A)$-measurable maps $(\gamma^{n,m})$, uniformly bounded in $n,m$.
Plug this decomposition of $F$ into \reff{barYnm}, and let us consider the process $\{\Gamma_{ts}^{n,m},t\leq s\leq T\}$ of dynamics:
\beqs
d\Gamma_{ts}^{n,m} &=& \Gamma_{ts}^{n,m}[ (\alpha_s^{n,m}-n) ds + \beta_s^{n,m} dW_s],  \;\;\; t \leq s \leq T, \;\;\; \Gamma_{tt}^{n,m} = 1,
\enqs
and given explicitly by:
\beqs
\Gamma_{ts}^{n,m} &=& e^{-n(s-t)} e^{\int_t^s \alpha_u^{n,m} du} M_{ts}^{n,m}, \;\;\;
M_{ts}^{n,m} \; = \; \frac{\Ec_s\big(\int_0^. \beta_u^{n,m} dW_u \big)}{\Ec_t\big(\int_0^. \beta_u^{n,m} dW_u \big)}, \; t \leq s \leq T,
\enqs
where $\Ec_t(\cdot)$ is the Dol\'eans-Dade exponential.  Since $\beta^{n,m}$ is a bounded process, we see that $\{M_{ts}^{n,m}$, $t\leq s\leq T\}$ is a uniformly integrable martingale, with
$M_{tT}^{n,m}$ $\in$ ${\bf L^p(}\Fc_T{\bf)}$  for any $p$ $\geq$ $1$.
%Thus, $\Gamma_{tT}^{n,m}$ also lies in ${\bf L^p(\Omega,\Fc_T)}$, for any $p$ $\geq$ $1$ since $\alpha^{n,m}$ is bounded.
%We can then rewrite $\bar Y^{n,m,k}$ as:
%\beqs
%\bar Y_t^{n,m,k}  &=&   \xi - U_T^k +  \int_t^T \big( F(s,0,0,0) + \alpha_s^{n,m} U_s^k + \beta_s^{n,m} v_s^k + u_s^k - \delta_s^{n,m} \big) ds  \\
%& & \; + \;  \int_t^T \big[ \alpha_s^{n,m} \bar Y_s^{n,m,k} - n (U_s - U_s^k -  \bar Y_s^{n,m,k})_- \big] ds \\
%& & \; + \; \int_t^T \int_A \big[ \gamma_s^{n,m}(a) L_s^{n,m}(a)  + m  (L_s^{n,m}(a))_+ - \nu_s(a) L_s^{n,m}(a) \big] \lambda(da) ds  \\
%& & \; - \;  \int_t^T  \big(Z_s^{n,m} -  v_s^k\big) (dW_s - \beta_s^{n,m} ds)  -  \int_t^T\int_A  L_s^{n,m}(a) \tilde\mu^\nu(ds,da).
%\enqs
By applying It\^o's formula to the product $\{\Gamma_{ts}^{n,m}\bar Y_s^{n,m,k},t\leq s\leq T\}$, we then obtain:
%$e^{\int_0^t(\alpha_s^{n,m} - n) ds} \bar Y_t^{n,m}$, we then obtain by setting $\Gamma_{ts}^{n,m}$ $=$ $e^{\int_t^s(\alpha_u^{n,m} - n) du}$:
\beqs
\bar Y_t^{n,m,k}  &=&  \Gamma_{tT}^{n,m} \big(  \xi - U_T^k \big) + \int_t^T \Gamma_{ts}^{n,m} \big( F(s,0,0,0) + \alpha_s^{n,m} U_s^k + \beta_s^{n,m} \var_s^k + \ups_s^k  \big) ds \\
& & \; + \;  \int_t^T  \Gamma_{ts}^{n,m} \big[ n \bar Y_s^{n,m,k} - n (U_s - U_s^k -  \bar Y_s^{n,m,k})_- - \delta_s^{n,m} \big] ds \\
& & \; + \; \int_t^T \int_A  \Gamma_{ts}^{n,m} \big[ \gamma_s^{n,m}(a) L_s^{n,m}(a)  + m  (L_s^{n,m}(a))_+ - \nu_s(a) L_s^{n,m}(a) \big] \lambda(da) ds  \\
& & \; - \;  \int_t^T  \Gamma_{ts}^{n,m}  \big(Z_s^{n,m} -  \var_s^k + \bar Y_s^{n,m,k} \beta_s^{n,m}\big) dW_s   -  \int_t^T\int_A  \Gamma_{ts}^{n,m}   L_s^{n,m}(a) \tilde\mu^\nu(ds,da),
\enqs
for any $\nu$ $\in$ $\Vc$, where we introduced the compensated measure $\tilde\mu^\nu$ of $\mu$ under $\P^\nu$.
By choosing $\nu$ $=$ $\nu^{n,m,\eps}$  $\in$ $\Vc$ defined by: $\nu_t^{n,m,\eps}(a)$ $=$ $(\gamma_t^{n,m}(a) +m)1_{\{L_t^{n,m}(a)\geq 0\}} + (\gamma_t^{n,m}(a)+\eps)1_{\{L_t^{n,m}(a)<0\}}$, for some arbitrary
$\eps$ $>$ $0$, we see that:
\beqs
\gamma_t^{n,m}(a) L_t^{n,m}(a)  + m  (L_t^{n,m}(a))_+ - \nu_t^{n,m}(a) L_t^{n,m}(a) &=& - \eps L_t^{n,m}(a) 1_{\{L_t^{n,m}(a)<0\}}.
%& \rightarrow & 0, \;\;\; a.s
\enqs
%when $\eps$ goes to zero.
Observe also  that
\beqs
n \bar Y_t^{n,m,k} - n (U_t - U_t^k -  \bar Y_t^{n,m,k})_-  - \delta_s^{n,m} & \leq & 0, \;\;\; 0 \leq t\leq T, \;  a.s.
\enqs
since $U$ $\leq$ $U^k$, and $\delta^{n,m}$ $\geq$ $0$. Recalling that $\xi\leq U_T\leq U_T^k$,  the explicit expression of $\Gamma^{n,m}$,  and  the fact that
$(\alpha^{n,m})$, $(\beta^{n,m})$ are uniformly bounded in $(t,\omega,n,m)$,  we then  get the existence of some positive constant $C$ such that:
\beq
\bar Y_t^{n,m,k}  & \leq &  C  \int_t^T e^{-n(s-t)}  M_{ts}^{n,m}  \big( |F(s,0,0,0)| + |U_s^k| + |\var_s^k| + |\ups_s^k|  \big) ds \label{inegYnm} \\
& & \!\!\! - \;  \eps \int_t^T \int_A \Gamma_{ts}^{n,m}  L_s^{n,m}(a) 1_{\{L_s^{n,m}(a)<0\}}  \lambda(da) ds  \nonumber \\
& & \!\!\! - \int_t^T  \Gamma_{ts}^{n,m}  \big(Z_s^{n,m} -  \var_s^k + \bar Y_s^{n,m,k} \beta_s^{n,m} \big) dW_s    -  \int_t^T\int_A \Gamma_{ts}^{n,m}   L_s^{n,m}(a) \tilde\mu^{\nu^{n,m,\eps}}(ds,da),
\nonumber
\enq
for any $n,m,k$ $\in$ $\N\setminus\{0\}$, $\eps$ $>$ $0$. Denote by $S_{t}^{n,m,k}$ $=$ $\int_0^t  \Gamma_{0s}^{n,m}  \big(Z_s^{n,m} -  \var_s^k + \bar Y_s^{n,m,k} \beta_s^{n,m} \big) dW_s$, $0\leq t\leq T$, which is
a $\P^\nu$-local martingale, for any $\nu$ $\in$ $\Vc$,  by recalling that $W$ remains a Brownian motion under $\P^\nu$.  From Burkholder-Davis-Gundy, Bayes formula,  Cauchy-Schwarz,
and Doob inequalities, we have
\beq
& & \E^\nu \big[ \sup_{0\leq t\leq T} |S_t^{n,m,k}| \big]  \nonumber  \\
& \leq & C \E^\nu \big[ \sqrt{< S^{n,m,k}>_T} \big]
 =  C \E^\nu \big[  \sqrt{\int_0^T |\Gamma_{0t}^{n,m}|^2  | Z_t^{n,m} -  \var_t^k + \bar Y_t^{n,m,k} \beta_t^{n,m}|^2 dt}  \big] \nonumber \\
& \leq & C \E \Big[ \zeta_T^\nu    \sup_{0\leq t\leq T}  \Gamma_{0t}^{n,m} \sqrt{\int_0^T  | Z_t^{n,m} -  \var_t^k + \bar Y_t^{n,m,k} \beta_t^{n,m}|^2  dt}  \Big]  \nonumber \\
& \leq & C \Big(\E \big[ |\zeta_T^\nu|^4 \big]  \E \big[  \sup_{0\leq t\leq T}  |\Gamma_{0t}^{n,m}|^4 \big] \Big)^{\frac{1}{4}}  \sqrt{ \E \big[ \int_0^T  | Z_t^{n,m} -  \var_t^k + \bar Y_t^{n,m,k} \beta_t^{n,m}|^2  dt \big] }  \nonumber \\
& \leq & C \Big(\E \big[ |\zeta_T^\nu|^4 \big]  \E \big[  |M_{0T}^{n,m}|^4 \big] \Big)^{\frac{1}{4}}  \sqrt{ \E \big[ \int_0^T  | Z_t^{n,m} -  \var_t^k + \bar Y_t^{n,m,k} \beta_t^{n,m}|^2  dt \big] }  \nonumber \\
& <& \infty,  \label{truemar}
\enq
where we used the fact that $\alpha^{n,m}$, $\beta^{n,m}$ are bounded processes, $Z^{n,m}$, $\vartheta^k$ lie in ${\bf L^2(W)}$, and $\bar Y^{n,m,k}$ in ${\bf L^2(0,T)}$.
Therefore, $S^{n,m,k}$ is a uniformly $\P^\nu$-integrable martingale for any $\nu$ $\in$ $\Vc$, and similarly we show that
$\int_0^t \int_A \Gamma_{ts}^{n,m}   L_s^{n,m}(a) \tilde\mu^{\nu}(ds,da)$ is a  $\P^\nu$-martingale. Hence, by taking conditional expectation with respect to $\P^{\nu^{n,m,\eps}}$ into \reff{inegYnm}, we
have for all $n,m,k$ $\in$ $\N\setminus\{0\}$, $\eps$ $>$ $0$:
\beq
\bar Y_t^{n,m,k}  & \leq &  \frac{C}{n}    \E^{\nu^{n,m,\eps}} \Big[ \sup_{t\leq s\leq T}  M_{ts}^{n,m}  \big( |F(s,0,0,0)| + |U_s^k| + |\var_s^k| + |\ups_s^k|  \big) \big| \Fc_t \Big] \nonumber \\
& & \; - \eps  \; \E^{\nu^{n,m,\eps}} \Big[  \int_t^T \int_A \Gamma_{ts}^{n,m}  L_s^{n,m}(a) 1_{\{L_s^{n,m}(a)<0\}}  \lambda(da) ds \big| \Fc_t  \Big]  \nonumber \\
&\leq &   \frac{C}{n}    \esssup_{\nu\in\Vc} \E^{\nu} \Big[ \sup_{t\leq s\leq T}  M_{ts}^{n,m}  \big( |F(s,0,0,0)| + |U_s^k| + |\var_s^k| + |\ups_s^k|  \big) \big| \Fc_t \Big] \label{Ynmeps}  \\
& & \; + \eps  \; \E \Big[ \frac{\zeta^{\nu^{n,m,\eps}}_T}{\zeta^{\nu^{n,m,\eps}}_t}  \int_t^T \int_A \Gamma_{ts}^{n,m}  |L_s^{n,m}(a)|   \lambda(da) ds \big| \Fc_t  \Big], \;\;\; 0 \leq t \leq T,
\nonumber
\enq
from Bayes formula.  Now,  for $\eps$ $\leq$ $m$, we see that $\nu^{n,m,\eps}$ $\leq$ $\bar\nu^{n,m}$ $:=$ $\gamma^{n,m} + m$, and so:
\beq \label{boundzeta}
0 \leq \; \frac{\zeta^{\nu^{n,m,\eps}}_T}{\zeta^{\nu^{n,m,\eps}}_t} & \leq & \frac{\zeta^{\bar\nu^{n,m}}_T}{\zeta^{\bar\nu^{n,m}}_t}  \exp\bigg(\int_t^T \int_A \bar\nu^{n,m}_s(a)\lambda(da)ds \bigg).
\enq
This shows that
\beq \label{Reps}
\lim_{\eps\rightarrow 0} \eps \;  \E \Big[ \frac{\zeta^{\nu^{n,m,\eps}}_T}{\zeta^{\nu^{n,m,\eps}}_t}  \int_t^T \int_A \Gamma_{ts}^{n,m}  |L_s^{n,m}(a)|   \lambda(da) ds \big| \Fc_t  \Big] &=& 0, \;\;\; 0 \leq t\leq T,
\enq
and so by sending $\eps$  to zero into \reff{Ynmeps}:
\beqs
& & (U_t^k - Y_t^{n,m})_- \; = \; (\bar Y_t^{n,m,k})_+ \\
& \leq & \frac{C}{n}   \esssup_{\nu\in\Vc} \E^{\nu} \Big[ \sup_{t\leq s\leq T}  M_{ts}^{n,m}  \big( |F(s,0,0,0)| + |U_s^k| + |\var_s^k| + |\ups_s^k|  \big) \big| \Fc_t \Big] \\
& \leq &  \frac{C}{n}  \esssup_{\nu\in\Vc} \E^{\nu} \Big[ \sup_{t\leq s\leq T}  |M_{ts}^{n,m}|^{\frac{p}{p-2}} +  \sup_{t\leq s\leq T} \big( |F(s,0,0,0)|^{\frac{p}{2}} + |U_s^k|^{\frac{p}{2}} + |\var_s^k|^{\frac{p}{2}}
+ |\ups_s^k|^{\frac{p}{2}}  \big) \big| \Fc_t \Big]
\enqs
for all $0\leq t\leq T$, and $p$ $>$ $2$, by Young  inequality.  Recall that $W$ is a Brownian motion under $\P^\nu$, and so
$\{M_{ts}^{n,m},t\leq s\leq T\}$ is a martingale under $\P^\nu$, for any $\nu$ $\in$ $\Vc$. By Doob's inequality, we then have with $q$ $=$ $p/(p-2)$ $>$ $1$:
\beqs
\E^{\nu} \Big[ \sup_{t\leq s\leq T}  |M_{ts}^{n,m}|^q \big| \Fc_t \Big] & \leq & \Big( \frac{q}{q-1}\Big)^q \E^\nu \big[ |M_{tT}^{n,m}|^q \big| \Fc_t \big]  \\
& \leq & \Big( \frac{q}{q-1}\Big)^q \exp\big( q(q-1) \|\beta\|^2_\infty (T-t) \big),
\enqs
where $\|\beta\|_\infty$ is a uniform bound of $(\beta^{n,m})$,  hence independent of $n,m$ and $\nu$ $\in$ $\Vc$. We then deduce that
\beqs
& & (U_t^k - Y_t^{n,m})_- \\
& \leq &  \frac{C}{n} \Big( 1 +  \esssup_{\nu\in\Vc} \E^{\nu}  \Big[ \sup_{t\leq s\leq T} \big( |F(s,0,0,0)|^{\frac{p}{2}} + |U_s^k|^{\frac{p}{2}} + |\var_s^k|^{\frac{p}{2}}
+ |\ups_s^k|^{\frac{p}{2}}  \big) \big| \Fc_t \Big]  \Big)
\enqs
for all $0\leq t\leq T$, $n,m,k$ $\in$ $\N\setminus\{0\}$. By Cauchy-Schwarz inequality,  we then obtain:
\beqs
& & \E \Big[ n \int_0^T (U_t^k - Y_t^{n,m})_- dt\Big]^2  \\
& \leq & C\bigg( 1 + \int_0^T \E \Big[ \esssup_{\nu\in\Vc} \E^{\nu}  \Big[ \sup_{t\leq s\leq T} \big( |F(s,0,0,0)|^{p} + |U_s^k|^{p} + |\var_s^k|^{p}
+ |\ups_s^k|^{p}  \big) \big| \Fc_t \Big] dt  \bigg).
\enqs
By taking   $p$ $>$ $2$ as in Assumption {\bf (H1)}(iii), and then sending $k$ to infinity in the l.h.s. of the above inequality, we get the required uniform estimate on
$K^{n,m,-}$.
 \ep

\vspace{2mm}

\begin{Corollary} \label{corollYm}
Let assumptions {\bf (H0)} and {\bf (H1)} hold.  Then, we have
\beqs
\sup_{m\in\N} \Big(   \|Y^{m}\| _{_{\bf S^2}} +   \|Z^{m}\|_{_{\bf L^2(W)}} + \|L^{m}\|_{_{\bf L^2(\tilde \mu)}} + \|K^{m,+}\|_{_{\bf S^2}} + \|K^{m,-}\|_{_{\bf S^2}} \Big)& < & \infty,
\enqs
where $K_t^{m,+}$ $:=$ $m\int_0^t\int_A \big(L_s^m(a)\big)_+ \lambda(da) ds$.
\end{Corollary}
{\bf Proof.}  From the bounds \reff{lowerYnm} and \reff{upperYnm}, we already have the uniform estimate for $\|Y^{m}\| _{_{\bf S^2}}$.
Moreover, by Lemmata  \ref{L:ZLbdd} and \ref{lemKnm-}, we have the uniform estimates:
\beqs
\sup_{n,m\in\N} \Big(    \|Z^{n,m}\|_{_{\bf L^2(W)}} + \|L^{n,m}\|_{_{\bf L^2(\tilde \mu)}} + \|K^{n,m,+}\|_{_{\bf S^2}} + \|K^{n,m,-}\|_{_{\bf S^2}} \Big)& < & \infty,
\enqs
We deduce that the weak limits  $(Z^m,L^m,K^{m,-})$ of $(Z^{m,n},L^{m,n},K^{n,m,-})$  when $n$ goes to infinity, are also uniformly bounded in
${\bf L^2(W)}\times{\bf L^2(\tilde \mu)}\times{\bf S^2}$.  From the strong convergence of $L^{n,m}$ to $L^m$ in ${\bf L^p(\tilde \mu)}$, $1\leq p<2$, we see by definition of $K^{n,m,+}$ and $K^{m,+}$
that $K_T^{n,m,+}$  converges strongly to $K_T^{m,+}$ in ${\bf L^p(}\Fc_T{\bf)}$, when $n$ goes to infinity.  Moreover, since $(K_T^{n,m,+})_n$
is uniformly bounded in ${\bf L^2(}\Fc_T{\bf)}$, it also converges weakly to $K_T^{m,+}$ in ${\bf L^2(}\Fc_T{\bf)}$.  It follows that $(K^{m,+})_m$ inherits from $(K^{n,m,+})_{n,m}$ the
uniform estimate  in ${\bf S^2}$.
\ep

\vspace{3mm}

We can now state the main result of this section as a consequence of the monotonic convergence theorem stated in Appendix B,
which extends to the Brownian-Poisson filtration framework the result  of  Peng and Xu \cite{pengxu05}.

\begin{Theorem} \label{ThmExistence}
Let assumptions {\bf (H0)} and {\bf (H1)}  hold. Then there exists a minimal solution $(Y,Z,L,K^+,K^-)$  $\in$
${\bf S^2}\times{\bf L^2(W)}\times{\bf L^2(\tilde \mu)}\times{\bf K^2}\times{\bf K^2}$ to the reflected BSDE with nonpositive jumps \reff{BSDEgen}-\reff{Lcons}-\reff{Ucons}-\reff{Sko}, where:
\begin{itemize}
\item[\textup{(i)}] $Y$ is the increasing limit of $(Y^m)_m$.
\item[\textup{(ii)}] $(Z,L)$ is the strong $($resp. weak$)$ limit of $(Z^m,L^m)_m$  in ${\bf L^p(W)}\times{\bf L^p(\tilde \mu)}$, with $p\in[1,2)$, $($resp. in ${\bf L^2(W)}\times{\bf L^2(\tilde \mu)}$$)$.
\item[\textup{(iii)}] $K^+_t$ is the weak limit of $(K_t^{m,+})_m$ in ${\bf L^2}(\Fc_t)$,  and $K^-_t$ is the strong limit of $(K_t^{m,-})_m$ in ${\bf L^2}(\Fc_t)$, for any $0 \leq t \leq T$.
\end{itemize}
\end{Theorem}
{\bf Proof.}  We already know that $(Y^m)_m$ is a nondecreasing sequence in ${\bf S^2}$, which converges to some $Y$, which satisfies $\underline{Y}\leq Y\leq \bar Y$
from  \reff{lowerYnm} and \reff{upperYnm}, and so lies in ${\bf S^2}$.  By Lemma \ref{lemKm}  and Corollary \ref{corollYm}, we then see that
the sequence $(Y^m,Z^m,L^m,K^{m,+},K^{m,-})_m$ solution to the BSDE \reff{BSDEmaxm} satisfies all the conditions of the monotonic limit Theorem  \ref{MonLimThm}.  This provides the existence
of $(Z, L, K^+, K^-)$ $\in$ ${\bf L^2(W)}\times{\bf L^2(\tilde \mu)}\times{\bf K^2}\times{\bf K^2}$ as in the assertions (ii) and (iii)  of Theorem \ref{ThmExistence}  such that the quintuple $( Y, Z, L, K^+, K^-)$ solves \reff{BSDEgen}.

From the strong convergence in ${\bf L^1(\tilde\mu)}$ of $(L^m)_m$  to $L$, and since $\lambda(A)$ $<$ $\infty$, we have
\beqs
\E \Big[ \int_0^T \int_A \big( L_t^m(a)\big)_+ \lambda(da) dt \Big] & \longrightarrow & \E \Big[ \int_0^T \int_A \big( L_t(a)\big)_+ \lambda(da) dt \Big],
\enqs
as $m$ goes to infinity. Moreover, since $K_T^{m,+}$ $=$ $m\int_0^T(L_t(a))_+\lambda(da)dt$ is bounded in $m$ in ${\bf L^2(}\Fc_T{\bf)}$, this implies that
\beqs
\E \Big[ \int_0^T \int_A \big( L_t(a)\big)_+ \lambda(da) dt \Big] &=& 0,
\enqs
which means that the constraint \reff{Lcons} is satisfied. The upper reflection \reff{Ucons} is obviously satisfied from \reff{consUm} and by sending $m$ to infinity.  Let us now  check the Skorohod reflecting condition \reff{Sko}. We recall  from \reff{Skom} that $\int_0^T (U_{t^-} - Y^m_{t^-})dK_t^{m,-}$ $=$ $0$.  Together  with the fact that
$U_{t^-} - Y^m_{t^-} \geq U_{t^-} - Y_{t^-} \geq 0$, this  yields $\int_0^T (U_{t^-} - Y_{t^-})dK_t^{m,-} = 0$. Since $(K^{m,-})_m$ converges strongly to $K^-$ in ${\bf S^2}$, this implies that the measure
$dK^{m,-}$ converges weakly to $dK^-$, and so  $\int_0^T (U_{t^-} - Y_{t^-})dK_t^{-} = 0$ a.s.

It remains to prove the minimality condition. Let $(\tilde Y, \tilde Z,\tilde L,\tilde K^+,\tilde K^-)$ be another solution to the reflected BSDE with nonpositive jumps \reff{BSDEgen}-\reff{Lcons}-\reff{Ucons}-\reff{Sko}.
We then see that  $\int_0^t \int_A (\tilde L_s(a))_+\lambda(da)ds$ $=$ $0$, and thus $F(t,\tilde Y_t,\tilde Z_t,\tilde L_t)$ $=$ $F_m(t,\tilde Y_t,\tilde Z_t,\tilde L_t)$, for $0\leq t\leq T$.
From  the comparison Theorem \ref{CompThm2}, we deduce that  $Y_t^m$ $\leq$ $\tilde Y_t$, $0 \leq t \leq T$. Taking the limit with respect to $m$, this proves the minimality condition:  $Y_t \leq \tilde Y_t$, $0 \leq t \leq T$.
\ep

\vspace{2mm}

\begin{Remark} \label{limitorder}
{\rm  The order of the limits: first let $n$ to infinity, and then let $m$ to infinity, is crucial in our approach. Indeed, by sending first $n$ to infinity, we
get a nondecreasing sequence of processes $(K^{m,-})_m$ (see Lemma \ref{lemKm}), which is a  required property  for applying the monotonic convergence theorem in Theorem \ref{ThmExistence}. On the other hand, if we would first let $m$ to infinity in the double sequence
$(Y^{n,m},Z^{n,m}, L^{n,m},K^{n,m,+},K^{n,m,-})$,  then we would obtain a minimal solution $(\hat Y^n,\hat Z^n,\hat K^{n,+})$ to the BSDE with nonpositive jumps:
\beq
\hat Y_t^n & = & \xi + \int_t^T F(s,\hat Y_s^n,\hat Z_s^n,\hat L_s^n) ds - n \int_t^T (U_s  - \hat Y_s^n)_- ds + \hat K_T^{n,+}-\hat K_t^{n,+}  \nonumber \\
& & \;\; - \; \int_t^T \hat Z_s^n dW_s -  \int_t^T \int_A \hat L_s^n(a) \mu(ds,da), \;\;\;  0 \leq t \leq T, \label{BSDEYn}  \\
\hat L_t^n(a) & \leq & 0,\;\;\;  d\P \otimes dt \otimes\lambda(da) \; a.e. \nonumber
\enq
and $(\hat Y^n)_n$ is a nonincreasing sequence, converging to some $\hat Y$ $\geq$ $Y$ by \reff{Ynmcroidecroi}.
But neither $K^{n,+}$, which is the weak limit of $K^{n,m,+}$, as $m$ goes to infinity, nor $K_t^{n,-}$ $:=$ $n \int_0^t (U_s  - \hat Y_s^n)_- ds$,
satisfy monotonicity properties in $n$, which prevents to apply the monotonic convergence theorem to the sequence
$(\hat Y^n,\hat Z^n,\hat K^{n,+},\hat K^{n,-})_n$, and thus to identify $\hat Y$ $=$ $Y$ as the minimal solution to the reflected BSDE with nonpositive jumps.
This differs from the case of doubly reflected BSDEs where one can send  indifferently first  $m$ or $n$
to infinity.
\ep
}
\end{Remark}

\section{Dual game representation}

\setcounter{equation}{0} \setcounter{Assumption}{0}
\setcounter{Theorem}{0} \setcounter{Proposition}{0}
\setcounter{Corollary}{0} \setcounter{Lemma}{0}
\setcounter{Definition}{0} \setcounter{Remark}{0}

In this section, we consider the case where the generator $F(t,\omega)$ does not depend on $y,z,\ell$, and we provide a dual game representation of the minimal solution to the reflected BSDE with nonpositive jumps in terms of a family of equivalent probability measures and discount factors.
In addition to the set of probability measures $\P^\nu$, $\nu$ $\in$ $\Vc$ $=$ $\cup_m \Vc_m$ defined in the previous section, let us introduce for any
$n$ $\in$ $\N$, the set $\Theta_n$ of  $\F$-progressively measurable processes valued in $[0,n]$, and set $\Theta$ $=$ $\cup_n \Theta_n$, which shall represent the set of discount processes.
Inspired by Proposition 6.2 in \cite{cvikar96} and the dual representation in Section 4 of \cite{khapha12}, we prove an explicit  representation formula for  the minimal solution to the reflected BSDE with nonpositive jumps.

\begin{Proposition} \label{P:Representation}
\textup{(i)} For any $n$ $\in$ $\N$ and $m$ $\in$ $\N\setminus\{0\}$, the solution to the penalized BSDE \reff{penBSDE} admits the following dual representation formula:
\beqs
Y_t^{n,m} \!\!\! &=& \!\!\! \esssup_{\nu\in\Vc_m} \essinf_{\theta\in\Theta_n} G_t(\nu,\theta)
\; = \;  \essinf_{\theta\in\Theta_n}  \esssup_{\nu\in\Vc_m}  G_t(\nu,\theta),
\enqs
for all $0 \leq t \leq T$, where
\beqs
G_t(\nu,\theta) & := & \E^\nu \Big[ e^{-\int_t^T \theta_s ds} \xi   + \; \int_t^T  e^{-\int_t^s \theta_r dr} \big(F(s) + \theta_s U_s\big) ds \;\big| \;  \Fc_t \Big]. \nonumber
\enqs
\textup{(ii)} Under assumptions {\bf (H0)} and {\bf (H1)}, the minimal solution to the reflected BSDE with nonpositive jumps \reff{BSDEgen}-\reff{Lcons}-\reff{Ucons}-\reff{Sko} is explicitly represented as:
\beq \label{dualY}
Y_t &=& \esssup_{\nu\in\Vc} \essinf_{\theta\in\Theta} G_t(\nu,\theta), \;\;\; 0 \leq t \leq T.
\enq
\end{Proposition}
{\bf Proof.} (i) Fix $n$ $\in$ $\N$ and $m$ $\in$ $\N\setminus\{0\}$. For  $\theta$ $\in$ $\Theta$,  by applying It\^o's rule to the product of the processes $e^{-\int_0^\cdot\theta_s ds}$ and
 $Y^{n,m}$ in \reff{penBSDE}, and by introducing the compensated measure $\tilde\mu^\nu(dt,da)$ under $\P^\nu$ for $\nu$ $\in$ $\Vc$,   we obtain:
\beqs
Y_t^{n,m} &=& e^{-\int_t^T\theta_s ds} \xi  + \int_t^T e^{-\int_t^s \theta_r dr} \big(F(s) + \theta_s U_s\big) ds   \\
& & + \int_t^T \int_A e^{-\int_t^s \theta_r dr}  \big( m (L_s^{n,m}(a))_+  -  \nu_s(a) L_s^{n,m}(a)  \big) \lambda(da) ds \\
& & -  \int_t^T e^{-\int_t^s \theta_r dr} \big( n(U_s - Y_s^{n,m} )_- +   \theta_s(U_s - Y_s^{n,m})   \big) ds \\
& & - \int_t^T e^{-\int_t^s \theta_r dr}Z_s^{n,m} dW_s  -   \int_t^T \int_A e^{-\int_t^s \theta_r dr} L_s^{n,m}(a)  \tilde\mu^{\nu}(ds,da).
\enqs
By same arguments as in \reff{truemar} (see also Lemma 4.2 in \cite{khapha12}), we can check that the  $\P^\nu$ local martingales
$\{\int_t^s e^{-\int_t^u \theta_r dr}Z_u^{n,m} dW_u, t\leq s\leq T\}$ and    $\{\int_t^s \int_A e^{-\int_t^u \theta_r dr} L_u^{n,m}(a)  \tilde\mu^{\nu}(du,da), t\leq s\leq T\}$ are actually uniformly integrable $\P^\nu$-martingales, so that
by taking conditional expectation under $\P^\nu$:
\beq
Y_t^{n,m}  &=&   G_t(\nu,\theta)  +   \E^\nu \Big[ \int_t^T \int_Ae^{-\int_t^s \theta_r dr}  \big( m(L_s^{n,m}(a))_+ -\nu_s(a) L_s^{n,m}(a)  \big) \lambda(da) ds \nonumber \\
& & \hspace{2cm} -  \;  \int_t^T e^{-\int_t^s \theta_r dr}  \big( n(U_s - Y_s^{n,m} )_- +   \theta_s(U_s - Y_s^{n,m})   \big)  ds \big| \Fc_t \Big], \label{P:Ynm}
\enq
and this relation holds  for any $\nu$ $\in$ $\Vc$, and $\theta$ $\in$ $\Theta$.
Now, observe that for any $\nu$ $\in$ $\Vc_m$, hence valued in $(0,m]$, we have
\beqs
m (L_t^{n,m}(a))_+ - \nu_t(a) L_t^{n,m}(a) & \geq & 0, \;\;\;  0 \leq t \leq T, \;  a \in A, \; a.s.
\enqs
and for $\nu$ $=$ $\nu^\eps$  $\in$ $\Vc_m$ defined by: $\nu^\eps_t(a)$ $=$ $m1_{\{L_t^{n,m}(a)\geq 0\}} + \eps 1_{\{L_t^{n,m}(a)< 0\}}  $, for arbitrary $\eps$ $\in$ $(0,m]$, we have
\beqs
m (L_t^{n,m}(a))_+ - \nu_t^\eps(a) L_t^{n,m}(a) & = & - \eps L_t^{n,m}(a)  1_{\{L_t^{n,m}(a)< 0\}} , \;\;\; 0 \leq t \leq T, \; a \in A, \; a.s.
\enqs
Similarly,  for any  $\theta$ $\in$ $\Theta_n$, hence valued in $[0,n]$, we have
\beqs
n(U_t - Y_t^{n,m} )_-  +  \theta_t(U_t-Y_t^{n,m}) & \geq & 0, \;\;\;  0 \leq t \leq T, \; a.s.
\enqs
and for $\theta^*$ $\in$ $\Theta_n$ defined by: $\theta^*_t$ $=$ $n 1_{\{Y_t^{n,m}\geq U_t\}}$, we have
\beqs
n(U_t-Y_t^{n,m})_-  +  \theta_t^*(U_t-Y_t^{n,m}) & = & 0, \;\;\; 0 \leq t \leq T, \; a.s.
\enqs
Therefore, by \reff{P:Ynm}, we get
\beq
G_t(\nu,\theta^*) \; \leq \; Y_t^{n,m} &= & G_t(\nu^\eps,\theta^*) +  \eps R_t^{n,m,\eps}(\theta^*), \;\;\;  \forall \nu \in \Vc_m,  \label{G1} \\
& \leq & G_t(\nu^\eps,\theta) +  \eps R_t^{n,m,\eps}(\theta), \nonumber \\
& \leq &  G_t(\nu^\eps,\theta) +  \eps R_t^{n,m,\eps}(0), \;\;\;\;\; \forall \theta \in \Theta_n, \label{G2}
\enq
for all  $\eps$ $\in$ $(0,m]$, where we set:
\beqs
R_t^{n,m,\eps}(\theta) &:=&   \E^{\nu^\eps} \Big[ \int_t^T \int_Ae^{-\int_t^s \theta_r dr} | L_s^{n,m}(a) |    \lambda(da) ds \big| \Fc_t \Big].
\enqs
For fixed $m$, and by viewing  the BSDE \reff{penBSDE} as a penalized  BSDE  in $n$ for the upper-reflected BSDE
with generator $F_m$ in \reff{defFm}, we have by standard arguments based on It\^o's lemma, uniform estimates in $n$ for $(Y^{n,m},Z^{n,m},L^{n,m})$ in ${\bf S^2}\times{\bf L^2(W)}\times{\bf L^2(\tilde\mu)}$
(see Theorem 4.2 in \cite{essaky}). Actually, these arguments show that for all $0\leq t\leq T$,  there exists some real-valued  $\Fc_t$-measurable random variable $C_t^m$ such that
\beq
\label{EstimateL}
\sup_{n\in\N} \E \Big[ \int_t^T \int_A |L_s^{n,m}(a)|^2 \lambda(da) ds | \Fc_t \Big]  & \leq & C_t^m.
\enq
Moreover, since $\nu^\eps$ $\leq$ $m$, we see as in \reff{boundzeta} that $\zeta_T^{\nu^\eps}/\zeta_t^{\nu^\eps}$ $\leq$ $e^{m(T-t)\lambda(A)}\zeta_T^m/\zeta_t^m$, where
$\zeta^m$ is the Radon-Nikodym density of $d\P^\nu/d\P$ for $\nu$ $=$ $m$. Thus, by Cauchy-Schwarz inequality,  there exists some real-valued  $\Fc_t$-measurable random variable $\tilde C_t^m$ such that
\beq \label{Rm}
\sup_{n \in \N} R_t^{n,m,\eps}(0) & \leq & \tilde C_t^m,
\enq
for all $\eps$ $\in$ $(0,m]$. Now, by \reff{G1}, we have:  $\essinf_{\theta\in\Theta_n}  \esssup_{\nu\in\Vc_m}  G_t(\nu,\theta)$ $\leq$ $Y_t^{n,m}$, and by \reff{G2}, we get:
\beqs
Y_t^{n,m} & \leq & \esssup_{\nu\in\Vc_m} \essinf_{\theta\in\Theta_n}  G_t(\nu,\theta)  + \eps R_t^{n,m,\eps}(0).
\enqs
By \reff{Rm}, we see in particular that  $\eps R_t^{n,m,\eps}(0)$ $\rightarrow$ $0$ a.s.  as $\eps$ goes to zero.  Since we always have
$\esssup_{\nu\in\Vc_m} \essinf_{\theta\in\Theta_n}  G_t(\nu,\theta)$ $\leq$  $\essinf_{\theta\in\Theta_n}  \esssup_{\nu\in\Vc_m}  G_t(\nu,\theta)$, this shows that
\beq
Y_t^{n,m} \; = \; \lim_{\eps\rightarrow 0} G_t(\nu^\eps,\theta^*)  &= & \esssup_{\nu\in\Vc_m} \essinf_{\theta\in\Theta_n}  G_t(\nu,\theta) \nonumber \\
&=& \essinf_{\theta\in\Theta_n}  \esssup_{\nu\in\Vc_m}  G_t(\nu,\theta), \label{Ynminfsup}
\enq
i.e. $(\nu^\eps,\theta^*)$ $\in\Vc_m\times\Theta_n$ is an $\eps$-saddle point for $G_t(\nu,\theta)$.

\vspace{1mm}

\noindent (ii) By sending $m$ to infinity into \eqref{Ynminfsup}, and recalling that $Y^m$ $=$ $\lim_n Y^{n,m}$, we get:
\beq \label{Yminter}
Y_t^m &=& \essinf_{\theta\in\Theta}  \esssup_{\nu\in\mathcal{V}_m}  G_t(\nu,\theta) \; \geq \; \esssup_{\nu\in\mathcal{V}_m} \essinf_{\theta\in\Theta}  G_t(\nu,\theta).
\enq
On the other hand, for  arbitrary $n_0\in\mathbb{N}$, we see that for any $\theta\in\Theta_{n_0}$ and any $n\geq n_0$:
\beqs
n(U_t - Y_t^{n,m} )_-  +  \theta_t(U_t-Y_t^{n,m}) & \geq & 0, \;\;\;  0 \leq t \leq T, \; a.s.,
\enqs
which implies, from \eqref{P:Ynm},
\beq
\label{E:YnmG}
Y_t^{n,m}  &\leq&   G_t(\nu,\theta)  \\
& & \; + \;   \mathbb{E}^\nu \Big[ \int_t^T \int_Ae^{-\int_t^s \theta_r dr}  \big( m(L_s^{n,m}(a))_+ -\nu_s(a) L_s^{n,m}(a)  \big) \lambda(da) ds \big| \mathcal{F}_t \Big],  \nonumber
\enq
for any $\nu$ $\in$ $\mathcal{V}$, $\theta$ $\in$ $\Theta_{n_0}$, and $n\geq n_0$. Now note that, since $L^{n,m}$ $\rightarrow$ $L^m$ strongly in ${\bf L^p(\tilde\mu)}$, $p\in[1,2)$, then, up to a subsequence, $L^{n,m}$ $\rightarrow$ $L^m$ $d\mathbb{P}\otimes dt\otimes\lambda(da)$ almost everywhere. Moreover, as already recalled in step (i) of the proof, we have uniform estimates in $n$ for $(L^{n,m})\in{\bf L^2(\tilde\mu)}$, namely, from \eqref{EstimateL} with $t=0$,
\beq
\label{EstimateL0}
\sup_{n\in\mathbb{N}} \mathbb{E} \Big[ \int_0^T \int_A |L_s^{n,m}(a)|^2 \lambda(da) ds \Big]  & \leq & C_0^m,
\enq
for some positive constant $C_0^m$. Then, sending $n$ to infinity in \eqref{E:YnmG} we obtain, from Lebesgue's dominated convergence theorem,
\beq
\label{E:YmG}
Y_t^m  &\leq&   G_t(\nu,\theta)  \\
& & \; + \;   \mathbb{E}^\nu \Big[ \int_t^T \int_Ae^{-\int_t^s \theta_r dr}  \big( m(L_s^m(a))_+ -\nu_s(a) L_s^m(a)  \big) \lambda(da) ds \big| \mathcal{F}_t \Big],  \nonumber
\enq
for any $\nu$ $\in$ $\mathcal{V}$, $\theta$ $\in$ $\Theta_{n_0}$. Since $\Theta$ $=$ $\cup_n \Theta_n$, from the arbitrariness of $n_0$ we conclude that \eqref{E:YmG} remains true for all $\theta\in\Theta$.
Take $\tilde\nu^\varepsilon$  $\in$ $\mathcal{V}_m$ defined by: $\tilde\nu^\varepsilon_t(a)$ $=$ $m1_{\{L_t^m(a)\geq 0\}} + \varepsilon 1_{\{L_t^m(a)< 0\}}$, for arbitrary $\varepsilon$ $\in$ $(0,m]$, so that
\beqs
m (L_t^m(a))_+ - \nu_t^\varepsilon(a) L_t^m(a) & = & - \varepsilon L_t^m(a)  1_{\{L_t^m(a)< 0\}} , \;\;\; 0 \leq t \leq T, \; a \in A, \; a.s.,
\enqs
and  thus by  \eqref{E:YmG}:
\beq
\label{G3}
Y_t^m  &\leq& G_t(\tilde\nu^\varepsilon,\theta) + \varepsilon \tilde R_t^{m,\varepsilon}(\theta) \;\, \leq \;\, G_t(\tilde\nu^\varepsilon,\theta) + \varepsilon \tilde R_t^{m,\varepsilon}(0), \;\;\;\;\; \forall \theta \in \Theta,
\enq
for all  $\varepsilon$ $\in$ $(0,m]$, where we set:
\beqs
\tilde R_t^{m,\varepsilon}(\theta) &:=&   \mathbb{E}^{\tilde\nu^\varepsilon} \Big[ \int_t^T \int_Ae^{-\int_t^s \theta_r dr} | L_s^m(a) |    \lambda(da) ds \big| \mathcal{F}_t \Big].
\enqs
Using again the uniform estimate \eqref{EstimateL0} and the fact that, up to a subsequence, $L^{n,m}$ $\rightarrow$ $L^{m}$ $d\mathbb{P}\otimes dt\otimes \lambda(da)$ a.e., we obtain, from \eqref{EstimateL} and Lebesgue's dominated convergence theorem,
\beqs
\label{EstimateLm}
\mathbb{E} \Big[ \int_t^T \int_A |L_s^m(a)|^2 \lambda(da) ds | \mathcal{F}_t \Big]  & \leq & C_t^m.
\enqs
Moreover, as in step (i) of the proof, since $\tilde\nu^\varepsilon$ $\leq$ $m$ we see that $\zeta_T^{\tilde\nu^\varepsilon}/\zeta_t^{\tilde\nu^\varepsilon}$ $\leq$ $e^{m(T-t)\lambda(A)}\zeta_T^m/\zeta_t^m$. Thus, by Cauchy-Schwarz inequality, it follows that, for all $\varepsilon$ $\in$ $(0,m]$,
\beqs \label{Rm2}
\tilde R_t^{m,\varepsilon}(0) & \leq & \tilde C_t^m,
\enqs
with the same real-valued  $\mathcal{F}_t$-measurable random variable $\tilde C_t^m$ as in \eqref{Rm}. Then, from \eqref{G3} we get
\beqs
Y_t^m & \leq & \esssup_{\nu\in\mathcal{V}_m} \essinf_{\theta\in\Theta}  G_t(\nu,\theta) + \varepsilon \tilde C_t^m,
\enqs
for all $\varepsilon$ $\in$ $(0,m]$. By sending $\varepsilon$ to zero, and combining with \eqref{Yminter}, we obtain:
\beq
Y_t^m &=&   \essinf_{\theta\in\Theta}  \esssup_{\nu\in\mathcal{V}_m}  G_t(\nu,\theta) \nonumber  \\
&=& \esssup_{\nu\in\mathcal{V}_m} \essinf_{\theta\in\Theta}  G_t(\nu,\theta).  \label{Ymsupinf}
\enq
Finally, by sending $m$ to infinity into \eqref{Ymsupinf}, we obtain the dual relation \eqref{dualY} for $Y$ $=$ $\lim_m Y^m$.
\ep

\begin{Remark}
{\rm  We don't know in general  if one can switch in \reff{dualY} the  essential infimum and supremum. Actually,  by considering $\hat Y^n$ $=$ $\lim_m Y^{n,m}$  the minimal solution to the
BSDE with nonnegative jumps \reff{BSDEYn}, one could show by similar arguments as in the second part (ii) of Proposition \ref{P:Representation} that:
\beqs
\hat Y_t^n &=& \essinf_{\theta\in\Theta_n}  \esssup_{\nu\in\Vc}  G_t(\nu,\theta) \; = \; \esssup_{\nu\in\Vc} \essinf_{\theta\in\Theta_n} G_t(\nu,\theta),
\enqs
so that $\hat Y$ $:=$ $\lim_n \hat Y^n$ satisfies:
\beqs
\hat Y_t &=& \essinf_{\theta\in\Theta}  \esssup_{\nu\in\Vc}  G_t(\nu,\theta).
\enqs
However, as pointed out in Remark \ref{limitorder}, we cannot conclude whether  $\hat Y_t$ is equal or strictly greater than  $Y_t$.
\ep
}
\end{Remark}

\section{Connection with HJB  Isaacs equation for controller-and-stopper games}

\label{S:HJB}

\setcounter{equation}{0} \setcounter{Assumption}{0}
\setcounter{Theorem}{0} \setcounter{Proposition}{0}
\setcounter{Corollary}{0} \setcounter{Lemma}{0}
\setcounter{Definition}{0} \setcounter{Remark}{0}

In this section, we show how the minimal solution to our class of reflected BSDEs with nonpositive jumps provides a probabilistic representation (hence  a Feynman-Kac formula)
to fully  nonlinear variational inequalities of Hamilton-Jacobi-Bellman (HJB) Isaacs  type  arising in a controller/stopper game, when considering a suitable Markovian framework.

\subsection{The Markovian framework}

We are given two measurable functions $b:\R^d\times \R^q \rightarrow \R^d$ and $\sigma:\R^d\times \R^q \rightarrow \R^{d\times d}$ and we introduce the forward Markov regime-switching process $(X,I)$ in $\R^d\times\R^q$ governed by:
\beq
dX_t &=&  b(X_t,I_t) dt + \sigma(X_t,I_t) dW_t \label{FSDEX} \\
dI_t &=& \int_A (a-I_{t^-}) \mu(dt,da). \label{FSDEI}
\enq
Therefore, the coefficients $b$ and $\sigma$, appearing in the dynamics of the diffusion process $X$, change according to the pure jump process $I$, which is associated to the Poisson random measure $\mu$ on $\R_+\times A$. We make the following standard assumption on the forward coefficients $b$ and $\sigma$:

\vspace{3mm}

\ni {\bf (HFC)}  \hspace{7mm}  There exists a constant $C$ such that
\beqs
|b(x,a)-b(x',a')| + |\sigma(x,a)-\sigma(x',a')| &\leq& C \big(|x-x'| + |a-a'|\big),
\enqs
for all $x,x'\in\R^d$ and $a,a'\in\R^q$.

\vspace{3mm}

\ni It is well-known that under hypothesis {\bf (HFC)} there exists a unique solution $(X^{t,x,a},I^{t,a})$ $=$ $(X_s^{t,x,a},I_s^{t,a})_{t \leq s \leq T}$ to \reff{FSDEX}-\reff{FSDEI} starting from $(x,a)\in\R^d\times\R^q$ at time $s=t\in[0,T]$. Furthermore, we have the standard estimates: for all $p \geq 2$, there exists some constant $C_p$ such that
\beq
\label{EstimateXI}
\E\Big[\sup_{t \leq s \leq T}\big(|X_s^{t,x,a}|^p + |I_s^{t,a}|^p \big)\Big] &\leq& C_p\big(1 + |x|^p + |a|^p\big),
\enq
for all $(t,x,a)\in[0,T]\times\R^d\times\R^q$.

\begin{Remark}
{\rm Notice  that the constant $C_p$ in \reff{EstimateXI} depends only on $p$, $T$, and the growth linear condition of $b,\sigma$ in {\bf (HFC)}.  Since the dynamics \reff{FSDEX} of $X$ is not changed by the change of probability measure $\P^\nu$, $\nu$ $\in$ $\Vc$ (recall that $W$ remains a Brownian motion under $\P^\nu$),  we then see that for all $p$ $\geq$ $2$:
\beqs
\E^\nu \Big[ \sup_{s\leq r\leq T} \big( |X_r^{t,x,a}|^p + |I_r^{t,a}|^p \big) | \Fc_s \Big] & \leq & C_p \big( 1 + |X_s^{t,x,a}|^p + |I_s^{t,a}|^p), \;\;\; t \leq s \leq T,
\enqs
for all $\nu$ $\in$ $\Vc$, and thus:
\beq \label{estimXnu}
\int_t^T \E\Big[ \esssup_{\nu\in\Vc} \E^\nu\big[\sup_{s \leq r \leq T}\big(|X_r^{t,x,a}|^p + |I_r^{t,a}|^p\big) \big| \Fc_s \big] \Big]ds &\leq&  C_p(1 + |x|^p + |a|^p),
\enq
for all $(t,x,a)\in[0,T]\times\R^d\times\R^q$.
\ep
}
\end{Remark}

\vspace{3mm}

Regarding the reflected BSDE with nonpositive jumps, the terminal condition, the generator function,  and the barrier are given respectively by some continuous functions $g:\R^d$ $\rightarrow$ $\R$,
$f:\R^d\times \R^q\times\R\times\R^d$ $\rightarrow$ $\R$,  and $u$ $:$ $[0, T]\times\R^d$ $\rightarrow$ $\R$. We make the following assumptions on the BSDE coefficients:

\vspace{3mm}

\ni {\bf (HBC)}
\begin{itemize}
\item[(i)] The functions $g$, $f(\cdot,\cdot,0,0)$ and $u$ satisfy a polynomial growth condition:
\beqs
\sup_{x\in\R^d,a\in\R^q} \frac{|f(x,a,0,0)|}{1 + |x|^h + |a|^h} + \sup_{t\in [0,T],x\in\R^d} \frac{|g(x) |+ |u(t,x)|}{1+|x|^h}  &<& \infty,
\enqs
for some $h$ $\geq$ $0$.
\item[(ii)] There exists some constant $C$ such that:
\beqs
|f(x,a,y,z) - f(x,a,y',z')| &\leq& C\big(|y-y'| + |z-z'|\big),
\enqs
for all $x\in\R^d$, $a\in\R^q$, $y,y'\in\R$, $z,z'\in\R^d$.
\item[(iii)] $u(T,x)$ $\geq$ $g(x)$, for all $x\in\R^d$, and there exists a  nonincreasing sequence of functions $(u^k)_k$  lying in  $C^{1,2}([0,T]\times\R^d)$, and converging pointwisely to $u$ such that the following polynomial growth condition holds
\beqs
\sup_{k \in \N} \sup_{t \in [0, T], x\in\R^d} \frac{\left|\frac{\partial u^k}{\partial t}(t,x)\right| + |D_x u^k(t,x)| + |D_x^2 u^k(t,x)|}{1 + |x|^h } &<& \infty,
\enqs
for some $h \geq 0$.

\end{itemize}

\vspace{3mm}

In this Markovian framework,  the reflected BSDE with nonpositive jumps \reff{BSDEgen}-\reff{Lcons}-\reff{Ucons}-\reff{Sko} takes the  form:
\beq
Y_t & = & g(X_T) + \int_t^T f(X_s,I_{s},Y_s,Z_s) ds    +  K_T^+ - K_t^+ - (K_T^- - K_t^-)   \label{BSDEgenM}  \\
& &  -  \int_t^T  Z_s dW_s  -  \int_t^T\int_A   L_s(a)  \mu(ds,da), \;\;\;   0 \leq t \leq T,  \; a.s. \nonumber
\enq
with
\beq
L_t(a)  & \leq & 0\;, \;\;\;\;\; d\P\otimes dt\otimes\lambda(da) \;\;  a.e. \label{LconsM}
\enq
and
\beq
Y_t & \leq & u(t,X_t)\;, \;\;\;\;\; 0 \leq t \leq T,\;a.s. \label{UconsM} \\
\int_0^T (u(t,X_t) - Y_{t^-}) dK^-_t & = & 0\;, \qquad\qquad\;\;\; a.s. \label{SkoM}
\enq

Notice that under {\bf (HFC)} and {\bf (HBC)} the terminal condition $\xi(\omega)$ $=$ $g(X_T(\omega))$, the generator $F(t,\omega,y,z,\ell)$ $=$ $f(X_t(\omega),I_{t^-}(\omega),y,z)$, and
the barrier $U_t(\omega) =u(t,X_t(\omega))$  clearly  satisfy the standing assumptions 1-4 in Section 2.  Let us now discuss about  conditions {\bf (H0)} and {\bf (H1)} in the two following remarks.

\vspace{2mm}

\begin{Remark} \label{remH0}
{\rm  Condition {\bf (H0)} is satisfied in our Markovian framework.  Actually, it is shown in Lemma 5.1 in \cite{khapha12} that under {\bf (HFC)} and {\bf (HBC)}(i), (ii), there exists  for any initial condition
$(t,x,a)\in[0,T]\times\R^d\times\R^q$,  a  solution $\{(\bar Y_s^{t,x,a},\bar Z_s^{t,x,a},\bar L_s^{t,x,a},\bar K_s^{t,x,a,+}),t \leq s \leq T\}$  to the BSDE with nonpositive jumps \reff{H01}-\reff{H02} when
$(X,I)=\{(X_s^{t,x,a},I_s^{t,a}),t \leq s \leq T\}$, with $\bar Y_s^{t,x,a}$ $=$ $\bar v(s,X_s^{t,x,a})$ for some deterministic function $\bar v$ on $[0,T]\times\R^d$ satisfying
the polynomial growth condition:
\beqs
\sup_{(t,x)\in[0,T]\times\R^d} \frac{|\bar v(t,x)|}{1 + |x|^{r}} &<& \infty
\enqs
for some $r$ $\geq$ $2$.  Such solution is constructed by It\^o's lemma from a smooth supersolution to
\beqs
- \Dt{\bar v} - \sup_{a \in A} [ \Lc^a \bar v + f(\cdot,a,\bar v,\sigma\trans(\cdot,a)D_x\bar v) ] & \geq & 0, \;\;\; \mbox{ on } [0,T)\times\R^d \\
\bar v(T,x) & \geq &  g(x), \;\;\;\;\; x \in \R^d,
\enqs
where
\beqs
\Lc^a \varphi &=& b(x,a). D_x \varphi + \frac{1}{2} {\rm  tr}(\sigma\sigma\trans(x,a)D_x^2 \varphi),
\enqs
which can be chosen equal to $\bar v(t,x)$ $=$ $\bar C e^{\rho(T-t)}(1+|x|^{r})$, with $r$ $=$ $\max(2,h)$, for $\bar C$ and $\rho$ positive large enough.
\ep
}
\end{Remark}

\begin{Remark}
\label{R:H1Markovian}
{\rm We also observe that assumption {\bf (H1)} is  satisfied in the present framework. More precisely, given an initial condition $(t,x,a)\in[0,T]\times\R^d\times\R^q$, let us consider the process
$U^k$, $k$ $\in$ $\N$,  defined by:
\beqs
U_s^k &:= &   u^k(s,X_s^{t,x,a}), \;\;\; t \leq s \leq T.
\enqs
By It\^o's formula, $U^k$ is in the form of condition {\bf (H1)}(ii), with
\beqs
\ups_s^k &=& \frac{\partial u^k}{\partial t}(s,X_s^{t,x,a}) + b(X_s^{t,x,a},I_s^{t,a}). D_x u^k(s,X_s^{t,x,a}) \\
& & + \; \frac{1}{2}{\rm tr}\big(\sigma\sigma\trans(X_s^{t,x,a},I_s^{t,a}) D_x^2 u^k(s,X_s^{t,x,a})\big), \\
\var_s^k &=& D_x u^k(s,X_s^{t,x,a})\trans \sigma(X_s^{t,x,a},I_s^{t,a}),
\enqs
for all $t \leq s \leq T$, a.s., and we clearly see  from \textbf{(HFC)}, \textbf{(HBC)}(iii), and \reff{EstimateXI} that
%$u^k$ $\in$ ${\bf L^2([0,T]}$ and $v^k$ $\in$ {\bf L^2(W)}$.
\beqs
\E\Big[\int_t^T |\ups_s^k|^2 ds\Big] + \E\Big[\int_t^T |\var_s^k|^2 ds\Big] &<&  \infty.
\enqs
Moreover, by using \reff{estimXnu}, and again from the polynomial  growth conditions on $b$, $\sigma$, $F$ and $u^k$ in \textbf{(HFC)}, \textbf{(HBC)},  there exists some $p$ $>$ $2$ such that
\beqs
\sup_{k \in \N}  \int_t^T \E\Big[  \esssup_{\nu\in\Vc} \E^\nu\big[\sup_{s \leq r \leq T}\big(|U_r^k|^p + |\ups_r^k|^p + |\var_r^k|^p\big) \big| \Fc_s \big] \Big]ds  & &  \nonumber \\
\; + \; \int_t^T\E\Big[\esssup_{\nu\in\Vc}\E^\nu\big[\sup_{s \leq r \leq T}\big|f(X_r^{t,x,a},I_r^{t,a},0,0)\big|^p \big| \Fc_s  \big]\big] ds  & \leq & C_p(1+|x|^p + |a|^p).
\enqs
for all $(t,x,a)\in[0,T]\times\R^d\times\R^q$.
\ep
}
\end{Remark}

\vspace{3mm}

From Theorem \ref{ThmExistence}, we get, for any initial condition $(t,x,a)\in[0,T]\times\R^d\times\R^q$,
the existence of a  minimal solution $\{( Y_s^{t,x,a}, Z_s^{t,x,a}, L_s^{t,x,a}, K_s^{t,x,a,+}, K_s^{t,x,a,-}),t \leq s \leq T\}$ to the Markovian reflected BSDE with nonpositive jumps \reff{BSDEgenM}-\reff{LconsM}-\reff{UconsM}-\reff{SkoM} when $(X,I)=\{(X_s^{t,x,a},I_s^{t,a}),t \leq s \leq T\}$. Moreover, as we shall see in the next paragraph, this minimal solution is written in this Markovian context as: $Y_s^{t,x,a}=v(s,X_s^{t,x,a},I_s^{t,a})$, where $v$ is a real-valued deterministic function defined on $[0,T]\times\R^d\times \R^q$ by
\beq \label{v}
v(t,x,a) &:=& Y_t^{t,x,a}, \qquad (t,x,a)\in[0,T]\times\R^d\times\R^q.
\enq
We aim at proving that  this function $v$    does not depend actually  on the argument $a$ in the interior of $A$,
and is connected to the fully nonlinear variational inequality of  HJB Isaacs type:
\beq
\max \Big[ - \Dt{v}  - \sup_{a\in A}\big(\Lc^{a}v   +  f(\cdot,a,v,\sigma\trans(\cdot,a) D_x v) \big);
 v - u \Big] \!\!&=&\!\! 0,  \; \mbox{ on } [0,T)\times\R^d  \label{HJB} \\
 v(T,x) \!\!&=&\!\! g(x), \;\;\; x \in \R^d. \label{termcondHJB}
\enq

\subsection{Viscosity property of the penalized BSDE}

Let us consider the Markovian penalized BSDE associated to \reff{BSDEgenM}-\reff{LconsM}-\reff{UconsM}-\reff{SkoM}
\beq
Y_t^{n,m} &=&   g(X_T) + \int_t^T f(X_s,I_{s},Y_s^{n,m},Z_s^{n,m}) ds \label{penBSDEM} \\
& & + \; m \int_t^T \int_A \big(L_s^{n,m}(a) \big)_+ \lambda(da) ds   - n  \int_t^T\big(u(s,X_s)-Y^{n,m}_s\big)_- ds \nonumber \\
& & \;  - \; \int_t^T  Z_s^{n,m} dW_s  -  \int_t^T\int_A   L_s^{n,m}(a)  \mu(ds,da), \;\;\; 0 \leq t \leq T,   \nonumber
\enq
and denote by $\{(Y_s^{n,m,t,x,a},Z_s^{n,m,t,x,a},L_s^{n,m,t,x,a}),t \leq s \leq T\}$ the unique solution to \reff{penBSDEM} when $(X,I)=\{(X_s^{t,x,a},I_s^{t,a}),t \leq s \leq T\}$ for any initial condition
$(t,x,a)\in[0,T]\times\R^d\times\R^q$. From the Markov property of the jump-diffusion process $(X,I)$, we recall from \cite{barbucpar97} that $Y_s^{n,m,t,x,a}$ $=$ $v^{n,m}(s,X_s^{t,x,a},I_s^{t,a})$, $t \leq s \leq T$, for some deterministic function $v^{n,m}$ defined on $[0,T]\times\R^d\times\R^q$ by
\beq \label{vnm}
v^{n,m}(t,x,a) &:=& Y_t^{n,m,t,x,a}, \qquad (t,x,a)\in[0,T]\times\R^d\times\R^q.
\enq
Next, for fixed $m$, let us  consider the limiting BSDE of \reff{penBSDEM} as $n$ goes to infinity, that is the reflected BSDE:
\beq
Y_t^m &=&   g(X_T) + \int_t^T f(X_s,I_{s},Y_s^m,Z_s^m) ds  + \; m \int_t^T \int_A \big(L_s^m(a)\big)_+ \lambda(da) ds \label{BSDEmM} \\
& & - \; (K_T^{m,-} - K_t^{m,-}) - \; \int_t^T  Z_s^m dW_s  -  \int_t^T\int_A   L_s^m(a)  \mu(ds,da), \qquad 0 \leq t \leq T,\,a.s.  \nonumber
\enq
and
\beq
Y_t^m &\leq& u(t,X_t), \qquad 0 \leq t \leq T,\,a.s. \label{consUmM} \\
\int_0^T (u(t,X_t) - Y_{t^-}^m) dK_t^{m,-} &=& 0, \qquad a.s. \label{SkomM}
\enq
and denote by $\{(Y_s^{m,t,x,a},Z_s^{m,t,x,a},L_s^{m,t,x,a},K_s^{m,t,x,a,+}),t \leq s \leq T\}$ the unique solution to \reff{BSDEmM}-\reff{consUmM}-\reff{SkomM} when $(X,I)=\{(X_s^{t,x,a},I_s^{t,a}),t \leq s \leq T\}$ for any initial condition $(t,x,a)\in[0,T]\times\R^d\times\R^q$. Since $Y^{n,m,t,x,a}$ converges to $Y^{m,t,x,a}$ as $n$ goes to infinity, we see from \reff{vnm}  that $Y^{m,t,x,a}$ may be written as $Y_s^{m,t,x,a}$ $=$  $v^m(s,X_s^{t,x,a},I_s^{t,a})$, $t \leq s \leq T$, where $v^m$ is the deterministic function defined on $[0,T]\times\R^d\times\R^q$ by:
\beq \label{vm}
v^m(t,x,a) & :=& \lim_{n\rightarrow\infty} v^{n,m}(t,x,a) \; = \;  Y_t^{m,t,x,a}, \;\;\;  (t,x,a)\in[0,T]\times\R^d\times\R^q.
\enq
From the convergence of   $Y^{m,t,x,a}$ to the minimal solution $Y^{t,x,a}$, when $m$ goes to infinity, as stated in Theorem \ref{ThmExistence}, we deduce that
$Y^{t,x,a}$ has indeed  the form $Y_s^{t,x,a}$ $=$ $v(s,X_s^{t,x,a},I_s^{t,a})$, with a deterministic function $v$ defined as the pointwise (nondecreasing) limit of $(v^m)_m$:
\beq \label{defv2}
v(t,x,a) &:=& \lim_{m\rightarrow\infty} v^m(t,x,a)\; = \; Y_t^{t,x,a}, \;\;\;  (t,x,a) \in [0,T]\times\R^d\times\R^q.
\enq

From the bounds \reff{lowerYnm}-\reff{upperYnm}, we have for all $m$ $\in$ $\N$:
$\underline v(t,x,a)$ $\leq$ $v^m(t,x,a)$ $\leq$ $\bar v(t,x)$, $(t,x,a)$ $\in$ $[0,T]\times\R^d\times\R^q$, where $\underline v$ $:=$ $v^0$ is associated to the reflected BSDE  $Y^m$ for $m$ $=$ $0$, and
$\bar v$ is the supersolution as defined in Remark \ref{remH0}. By the  polynomial growth condition on $\bar v$, and also on $\underline v$ (see e.g.  Lemma 3.2 in \cite{dqs13}), we deduce that $v^m$, and thus also
$v$ by passing to the limit, satisfy a polynomial growth condition:  there exist some positive constant $C$ and some $p \geq 2$, such that, for all $m$ $\in$ $\N$:
\beq \label{growthcondvmv}
 |v^m(t,x,a)| + |v(t,x,a)| \leq C(1 + |x|^p + |a|^p),
\enq
for all $(t,x,a)\in[0,T]\times\R^d\times\R^q$. As expected, for fixed $m$,  the function $v^m$ associated to the reflected BSDE with jumps \reff{BSDEmM}-\reff{consUmM}-\reff{SkomM}
is connected to the integro-differential variational inequality:
\beq  \label{IPDEvm}
\max \Big[  - \Dt{v^m}(t,x,a) - \Lc^{a} v^m(t,x,a)  - f(x,a,v^m(t,x,a),\sigma\trans D_xv^m(t,x,a)) & & \\
\;\;\;\; - \; m \int_A  \big( v^{m}(t,x,a')-v^{m}(t,x,a) \big)_+  \lambda(da') \; ;  & &  \nonumber \\
v^m(t,x,a)-u(t,x) \big] &=& 0, \nonumber
\enq
for $(t,x,a)\in[0,T)\times\R^d\times\R^q$, together with the terminal condition:
\beq \label{termcondvm}
v^m(T,x,a) &=& g(x),  \;\;\;\;\; (x,a)\in\R^d\times\R^q.
\enq
More precisely, we have the following result, which may be proved by extending to the multidimensional case Lemma 3.1 and Theorem 3.4 of \cite{dqs13}, and by
using Theorem \ref{CompThm} as comparison theorem for BSDEs with jumps.

\begin{Proposition}
\label{P:ViscPropvm}
Let assumptions {\bf (HFC)} and {\bf (HBC)} hold. The function $v^m$ in \reff{vm} is a continuous viscosity solution to \reff{IPDEvm}-\reff{termcondvm}, i.e., it is continuous on $[0,T]\times\R^d\times\R^q$, a viscosity supersolution
$($resp. subsolution$)$ to \reff{termcondvm}, i.e.
\beqs
v^m(T,x,a) &\geq& (resp.\text{ }\leq) \;\, g(x,a)
\enqs
for any $(x,a)\in\R^d\times\R^q$, and a viscosity supersolution $($resp.\text{ }subsolution$)$ to \reff{IPDEvm}, i.e.
\beq
\label{viscpropvm}
\max \Big[- \Dt{\varphi}(t,x,a) - \Lc^{a} \varphi(t,x,a) - f(x,a, v^m(t,x,a),\sigma\trans(x,a) D_x \varphi(t,x,a)) & & \\
- m \int_A  \big(\varphi(t,x,a')-\varphi(t,x,a) \big)_+  \lambda(da'); & & \nonumber \\
v^m(t,x,a)-u(t,x) \big] \geq \quad (\text{resp. }\leq) &0& \nonumber
\enq
for any $(t,x,a)\in[0,T)\times\R^d\times\R^q$ and any $\varphi\in C^{1,2}([0,T]\times(\R^d\times\R^q))$ such that
\beq
\label{minvm}
(v^m-\varphi)(t,x,a) &=& \min_{[0,T]\times\R^d\times\R^q}(v^m-\varphi) \quad (resp.\text{ }\max_{[0,T]\times\R^d\times\R^q}(v^m-\varphi)).
\enq
\end{Proposition}

%\vspace{1mm}

\begin{Remark}
\label{R:v^m<u}
{\rm Notice  that
\beq
\label{E:v^m<u}
v^m(t,x,a) & \leq & u(t,x), \qquad \text{for all }(t,x,a)\in[0,T]\times\R^d\times\R^q.
\enq
Indeed, for any $(t,x,a)\in[0,T]\times\R^d\times\R^q$, since $Y_s^{m,t,x,a}$ $=$  $v^m(s,X_s^{t,x,a},I_s^{t,a})$, $t \leq s \leq T$, we deduce, from \reff{consUmM} that
\beqs
\E\bigg[\frac{1}{s-t}\int_t^s\big(v^m(r,X_r^{t,x,a},I_r^{t,a})-u(r,X_r^{t,x,a})\big)dr\bigg]  &\leq & 0
\enqs
for all $t < s \leq T$. Since $(X^{t,x,a},I^{t,a})$ is c\`{a}dl\`{a}g, in particular it is right-continuous at time $t$. Therefore, \reff{E:v^m<u}  follows from the continuity of $v^m$ and $u$.
\ep
}
\end{Remark}

\subsection{HJB Isaacs equation}

This paragraph is devoted to the derivation of the equation satisfied in the viscosity sense by the function $v$ in \reff{defv2}, by passing to the limit, as $m$ goes to infinity, in the equation satisfied by $v^m$.
The first step is to prove that $v$ does not depend on $a$, which is basically  a consequence of the nonpositive jump constraint:
\beqs
L_s^{t,x,a}(a')  &=& v(s,X_s^{t,x,a},a') - v(s,X_s^{t,x,a},I_{s^-}^{t,x,a}) \; \leq \; 0, \;\;\;\;\;    d\P\otimes ds\otimes\lambda(da') \; a.e.
\enqs
providing that the function $v$ is continuous. However, as we do not know a priori that the function $v$ is continuous, we shall rely on (discontinuous) viscosity solutions arguments as in \cite{khapha12},
and make the following conditions on the set $A$ and the intensity measure $\lambda$:

\vspace{3mm}

\ni {\bf (H$A$)}  \hspace{7mm}  The interior set $\mathring{A}$ of $A$ is connex, and $A$ $=$ Adh$(\mathring{A})$, the closure of its interior.

\vspace{3mm}

\ni {\bf (H$\lambda$)}
\begin{itemize}
\item[(i)] The measure $\lambda$ supports the whole set $\mathring{A}$: for any $a\in\mathring{A}$ and any open neighborhood $\Oc$ of $a$ in $\R^q$ we have $\lambda(\Oc\cap\mathring{A})>0$.
\item[(ii)] The boundary of $A$: $\partial A$ $=$ $A\backslash\mathring{A}$, is negligible with respect to $\lambda$, i.e., $\lambda(\partial A)$ $=$ $0$.
\end{itemize}

\vspace{1mm}

\begin{Proposition} \label{propva}
Let assumptions {\bf (HFC)}, {\bf (HBC)}, {\bf (H$A$)}, and {\bf (H$\lambda$)} hold. Then the function $v$ does not depend on the variable $a$ on $[0,T)\times\R^d\times\mathring A$:
\beq \label{va}
v(t,x,a) &=& v(t,x,a'), \;\;\;  a,a' \in \mathring A,
\enq
for all $(t,x)$ $\in$ $[0,T)\times\R^d$.
\end{Proposition}
{\bf Proof.} The proof borrows most arguments from section 5.3 in \cite{khapha12}, and we only report here the main steps and the points to be modified. First, we see from \reff{E:v^m<u}, and sending $m$ to infinity that:
\beq \label{vlequ}
v & \leq & u \;\;\; \mbox{ on } \; [0,T]\times\R^d\times\R^q.
\enq
We next show that  the function $v$ is a viscosity supersolution to:
\beq \label{viscoa}
- | D_a v(t,x,a) | &=& 0, \;\;\; (t,x,a) \in [0,T)\times\R^d\times\mathring A,
\enq
i.e., for any $(t,x,a)\in[0,T)\times\R^d\times\mathring{A}$ and any function $\varphi\in C^{1,2}([0,T]\times(\R^d\times\R^q))$ such that $(v-\varphi)(t,x,a)=\min_{[0,T]\times\R^d\times\R^q}(v-\varphi)$, we have
\beqs
-\big|D_a\varphi(t,x,a)\big| &\geq& 0, \;\;\; \text{i.e.}  \;\;\; D_a \varphi(t,x,a) \; = \; 0.
\enqs
Indeed, let $(t,x,a)\in[0,T)\times\R^d\times\mathring{A}$ and $\varphi\in C^{1,2}([0,T]\times(\R^d\times\R^q))$ such that $(v-\varphi)(t,x,a)=\min_{[0,T]\times\R^d\times\R^q}(v-\varphi)$. We may assume, without loss of generality, that $v(t,x,a)$ $=$ $\varphi(t,x,a)$, $(t,x,a)$ is a strict minimum point, and we distinguish two cases: (i)  $v(t,x,a)\geq u(t,x)$.  From \reff{vlequ}, we have
\beqs
\varphi(t,x,a') &\leq& v(t,x,a') \;\,\leq\;\, u(t,x), \qquad \forall \; a'\in \R^q
\enqs
and $\varphi(t,x,a)$ $=$ $v(t,x,a)$ $=$ $u(t,x)$. It follows that  $\varphi(t,x,a)$ $=$ $\max_{a'\in\R^q}\varphi(t,x,a')$,  which yields:  $D_a\varphi(t,x,a)$ $=$ $0$, since $a$ $\in$ $\mathring A$.
(ii)  $v(t,x,a)<u(t,x)$. Then, for $m$ large enough, we also  have $v^m(t,x,a)$ $<$ $u(t,x)$, and $(t,x,a)$ is a local minimum point of $v^m-\varphi$.
By the viscosity supersolution property of $v^m$ to \reff{viscpropvm},  this implies:
\beqs
- \Dt{\varphi}(t,x,a) - \Lc^{a} \varphi(t,x,a) - f(x,a, v^m(t,x,a),\sigma\trans(x,a) D_x \varphi(t,x,a)) & & \\
- m \int_A  \big(\varphi(t,x,a')-\varphi(t,x,a) \big)_+  \lambda(da') & \geq & 0.
\enqs
By sending $m$ to infinity, we conclude as in the proof of Lemma 5.3 in \cite{khapha12} that: $\int_A \big(\varphi(t,x,a')$ $-$ $\varphi(t,x,a) \big)_+  \lambda(da')$ $=$ $0$, which means under {\bf (H$\lambda$)} that
$\varphi(t,x,a)$ $=$ $\max_{a'\in\R^q}\varphi(t,x,a')$, i.e., $D_a\varphi(t,x,a)$ $=$ $0$.

Finally, by arguing exactly as in Lemma 5.4 and Proposition 5.2 of \cite{khapha12}, we obtain under the additional condition {\bf (H$A$)}  the non dependence of $v$ on $a$ $\in$ $\mathring A$ from the viscosity supersolution property to  \reff{viscoa}.
\ep

\vspace{3mm}

From Proposition \ref{propva}, we can define by misuse of notation the function $v$ on $[0,T)\times\R^d$ by:
\beqs
v(t,x) &=& v(t,x,a), \;\;\; (t,x) \in [0,T)\times\R^d,
\enqs
for any $a$ $\in$ $\mathring A$, and we see that $v$ satisfies a polynomial growth condition when $x$ goes to infinity by \reff{growthcondvmv}. We finally state the viscosity property of $v$ to the HJB Isaacs type equation \reff{HJB}-\reff{termcondHJB}. Recall the definition of lower semicontinuous envelope $v_*$, and upper semicontinuous envelope $v^*$:
\beqs
v_*(t,x)&=&\liminf_{\substack{(t',x')\rightarrow(t,x)\\t'<T}}v(t',x') \qquad \text{and} \qquad v^*(t,x)\;\,=\;\,\limsup_{\substack{(t',x')\rightarrow(t,x)\\t'<T}}v(t',x'),
\enqs
for all $(t,x)\in[0,T]\times\R^d$.

\begin{Theorem} \label{thmvisco}
Let assumptions {\bf (HFC)}, {\bf (HBC)}, {\bf (H$A$)}, and {\bf (H$\lambda$)} hold. Then $v$ is a viscosity solution to \reff{HJB}-\reff{termcondHJB} in the sense that it verifies:

\noindent (i) {\it Viscosity supersolution property}:
\beq \label{termvsur}
v_*(T,x) &\geq& g(x),
\enq
for any $x\in\R^d$, and
\beq
\max \Big[ - \Dt{\varphi}(t,x) - \sup_{a\in A}\Big(\Lc^{a}\varphi(t,x)  +  f\big(x,a,v_*(t,x),\sigma\trans(x,a) D_x \varphi(t,x)\big)\Big); & &\label{vsur}\\
v_*(t,x) - u(t,x) \big] & \geq  &0 \nonumber
\enq
for any $(t,x)\in[0,T)\times\R^d$ and any $\varphi\in C^{1,2}([0,T]\times\R^d)$ such that  $(v_*-\varphi)(t,x)$ $=$ $\min_{[0,T]\times\R^d}(v_*-\varphi)$

\vspace{1mm}

\noindent (ii) {\it Viscosity subsolution property}:
\beq \label{termvsous}
v^*(T,x) &\leq& g(x),
\enq
for any $x\in\R^d$, and
\beq
\max \Big[ - \Dt{\varphi}(t,x) - \sup_{a\in A}\Big(\Lc^{a}\varphi(t,x)  +  f\big(x,a,v^*(t,x),\sigma\trans(x,a) D_x \varphi(t,x)\big)\Big); & & \label{vsous} \\
v^*(t,x) - u(t,x) \big] & \leq  &0 \nonumber
\enq
for any $(t,x)\in[0,T)\times\R^d$ and any $\varphi\in C^{1,2}([0,T]\times\R^d)$ such that  $(v^*-\varphi)(t,x)$ $=$ $\max_{[0,T]\times\R^d}(v^*-\varphi)$.
\end{Theorem}
{\bf Proof.}  The proof is quite  similar to the proof detailed  in Section 5.4 of \cite{khapha12}, and  we report only the main arguments and  the points to be modified with respect to the proof in \cite{khapha12}.

\noindent $\bullet$  {\it Viscosity supersolution property \reff{vsur}}:   Since $v$ is the pointwise limit of the nondecreasing sequence  of continuous functions $(v^m)$,  and recalling \reff{va},
we know (see e.g. \cite{bar94}) that $v$ is lower semicontinuous and so:
\beqs
v(t,x) \; = \; v_*(t,x) & = & \lim_{m\rightarrow\infty} v^m(t,x,a), \;\;\;\;\;  \forall (t,x,a) \in  [0,T]\times\R^d\times\mathring A.  \label{vinf1}
\enqs
%for all $(t,x,a)$ $\in$ $[0,T]\times\R^d\times\mathring A$
%where  $\displaystyle{\liminf_{m\rightarrow\infty}}  {}_*v_m(t,x,a)$ $:=$  $\displaystyle{\liminf_{\tiny{\begin{array}{c}m\rightarrow\infty \\ (t',x',a')\rightarrow (t,x,a)\\
%'<T\end{array}}} }v^m(t',x',a')$.
Fix now  $(t,x)\in[0,T)\times\R^d$, and let $\varphi\in C^{1,2}([0,T]\times\R^d)$ such that  $(v_*-\varphi)(t,x)$ $=$ $\min_{[0,T]\times\R^d}(v_*-\varphi)$, with a strict minimum without loss of generality.
We already know from \reff{vlequ} that $v_*$ $\leq$ $u$, and so distinguish two cases: (1) if $v^*(t,x)$ $=$ $u(t,x)$, then the viscosity supersolution property of $v$ at $(t,x)$ is obviously satisfied.  (2) Otherwise,  if
$v_*(t,x)$ $<$ $u(t,x)$, then for any arbitrary fixed $a$ $\in$ $\mathring A$, we have  $v^m(t,x,a)$ $<$ $u(t,x)$, and $(t,x,a)$ is a local minimum point of $v^m-\varphi$, for $m$ large enough.
From the viscosity supersolution property  \reff{viscpropvm}  of $v^m$ at $(t,x,a)$ with the test function $\varphi$, we  then get:
\beqs
- \Dt{\varphi}(t,x) - \Lc^{a} \varphi(t,x) - f(x,a, v^m(t,x,a),\sigma\trans(x,a) D_x \varphi(t,x)) & \geq & 0.
\enqs
By sending $m$ to infinity, and since $a$ is arbitrary in $\mathring A$, together with the continuity of the coefficients $b$, $\sigma$, and $f$ in the variable $a$, we obtain the required viscosity supersolution inequality:
\beqs
- \Dt{\varphi}(t,x) - \sup_{a \in A} \Big( \Lc^{a} \varphi(t,x) + f(x,a, v_*(t,x),\sigma\trans(x,a) D_x \varphi(t,x)) \Big) & \geq & 0.
\enqs

\vspace{1mm}

\noindent $\bullet$  {\it Viscosity subsolution property \reff{vsous}}:  By \reff{vlequ}, we have: $v^*$ $\leq$ $u$ on $[0,T)\times\R^d$, and so  it remains to show the viscosity subsolution property of $v$ to:
\beqs
- \Dt{v}  - \sup_{a\in A}\Big(\Lc^{a} v(t,x)  +  f\big(x,a,v(t,x),\sigma\trans(x,a) D_x v(t,x)\big)\Big) & \leq & 0.
\enqs
This follows  by same arguments as in \cite{khapha12}  from the viscosity subsolution property of $v^m$  to:
\beqs
- \Dt{v^m}(t,x,a) - \Lc^{a} v^m(t,x,a) - f(x,a, v^m(t,x,a),\sigma\trans(x,a) D_x v^m(t,x,a)) & & \\
- m \int_A  \big(v^m(t,x,a')-v^m(t,x,a) \big)_+  \lambda(da') & \leq & 0,
\enqs
and by sending $m$ to infinity under {\bf (H$\lambda$)}(ii).

\vspace{1mm}

\noindent $\bullet$ Finally, the viscosity supersolution and subsolution inequalities \reff{termvsur}, \reff{termvsous} are proved by same arguments as in  \cite{khapha12}.
\ep

\vspace{3mm}

\begin{Remark}
{\rm {\it Zero-sum controller/stopper game}

\noindent  Let us consider the particular and important case where the generator $f(x,a)$ does not depend on $(y,z)$, and $u(t,x)$ $=$ $g(x)$.
In this case, the  nonlinear variational inequality \reff{HJB}-\reff{termcondHJB} is the
HJB Isaacs equation associated to the following zero-sum controller-and-stopper game: let us introduce  the controlled diffusion process in $\R^d$
\beq \label{Xcontrol}
dX_s^\alpha &=& b(X_s^\alpha,\alpha_s) ds + \sigma(X_s^\alpha,\alpha_s) dW_s,
\enq
where the control $\alpha$ $\in$ $\Ac$ is an $\F^W$-progressively measurable process, valued in $A$,  affecting  both drift and diffusion coefficient, possibly degenerate.  Here  $\F^W$ denotes  the natural filtration generated by the Brownian motion $W$. Notice that the laws  $\P^\alpha$ of $X^\alpha$ under $\P$, for $\alpha$ varying in $\Ac$,  belong to a non dominated set of probability measures.
Given $(t,x)$ $\in$ $[0,T]\times\R^d$, and $\alpha$ $\in$ $\Ac$, we denote by $\{X_s^{t,x,\alpha},t\leq s\leq T\}$ the solution to \reff{Xcontrol} starting from $x$ at $s$ $=$ $t$. Let us also define
$\Tc_{t,T}$ as  the set of all $\F^W$-stopping times valued in $[t,T]$ for $0\leq t\leq T$, and consider  $\Pi_{t,T}$ the set of  stopping strategies $\pi$ $:$ $\Ac$ $\mapsto$ $\Tc_{t,T}$ satisfying a non-anticipative condition
as defined in \cite{bayhua13}.  The upper and lower value functions of the controller/stopper game are given by:
\beqs
\overline V(t,x) &:=& \inf_{\pi\in\Pi_{t,T}} \sup_{\alpha\in\Ac} \E \Big[ \int_t^{\pi[\alpha]} f(X_s^{t,x,\alpha},\alpha_s) ds + g(X_{\pi[\alpha]}^{t,x,\alpha}) \Big], \\
\underline V(t,x) &:=&  \sup_{\alpha\in\Ac} \inf_{\tau\in\Tc_{t,T}} \E \Big[ \int_t^\tau f(X_s^{t,x,\alpha},\alpha_s) ds  + g(X_\tau^{t,x,\alpha}) \Big], \;\;\;\;\;  (t,x) \in [0,T]\times\R^d.
\enqs
It is shown in \cite{bayhua13} that this game has a value, i.e., $\overline V$ $=$ $\underline V$ $=$ $V$, and that $V$ is the unique viscosity solution  to \reff{HJB}-\reff{termcondHJB} satisfying a polynomial growth condition.
By combining this result with Theorem \ref{thmvisco}, this shows that  $v$ $=$ $V$. In other words, we have provided a representation of HJB Isaacs equation, arising in zero-sum controller/stopper game, including control on possibly degenerate diffusion coefficient,  in terms of minimal solution to reflected BSDE with nonpositive jumps.  Furthermore, by combining with the dual game representation in Proposition \ref{P:Representation},
we obtain an original representation for the value function of the controller-and-stopper game:
\begin{align*}
&\;\,\inf_{\pi\in\Pi_{0,T}} \sup_{\alpha\in\Ac} \E \Big[ \int_0^{\pi[\alpha]} f(X_t^{\alpha},\alpha_t) dt + g(X_{\pi[\alpha]}^{\alpha}) \Big] \; = \;
\sup_{\alpha\in\Ac} \inf_{\tau\in\Tc_{0,T}} \E \Big[ \int_0^\tau f(X_t^{\alpha},\alpha_t) dt  + g(X_\tau^{\alpha}) \Big] \\
=&\;\,\sup_{\nu \in \Vc} \inf_{\theta\in\Theta} \E^\nu \Big[ \int_0^T  e^{-\int_0^t \theta_s ds} \big(f(X_t,I_t) + \theta_t g(X_t) \big) dt +  e^{-\int_0^T \theta_t dt} g(X_T) \Big].
\end{align*}
\ep
}
\end{Remark}

\section{Conclusion}

We introduced in this paper a class of reflected BSDEs with nonpositive jumps and upper obstacle, and showed in the Markov case  its connection with fully nonlinear variational inequalities arising typically in controller-and-stopper games with control both on drift and diffusion term.
Such  representation suggests an original  approach for  probabilistic numeri\-cal schemes  of  HJB Isaacs equations by discretization and simulation of this reflected BSDE with nonpositive jumps.
From a theoretical point of view, an open problem is  to relate this class of BSDEs to general controller-and-stopper games in the non Markovian case.
A variation of our class of BSDEs would be to consider reflected BSDEs with nonpositive jumps and lower obstacle, which is related to $\sup\sup$ problem over control and stopping time, and in other words to
optimal stopping under nonlinear expectation. Actually, the proof  of existence of a minimal solution by a double penalization approach is simpler since it would involve the sum (instead of the difference) of two nondecreasing processes. Another possible extension is  the class of doubly reflected BSDEs with nonpositive jumps  motivated by Dynkin games under nonlinear expectation (see \cite{matetal13}).

\appendix

\setcounter{equation}{0} \setcounter{Assumption}{0}
\setcounter{Theorem}{0} \setcounter{Proposition}{0}
\setcounter{Corollary}{0} \setcounter{Lemma}{0}
\setcounter{Definition}{0} \setcounter{Remark}{0}

\renewcommand\thesection{Appendices}

\section{}

\renewcommand\thesection{\Alph{subsection}}

\renewcommand\thesubsection{\Alph{subsection}.}

\subsection{Comparison theorems for sub and supersolutions to BSDEs with jumps}

We provide in this section two comparison theorems for BSDEs with jumps. We first recall a comparison theorem for sub and supersolutions to BSDEs driven by the Brownian motion $W$ and the Poisson random measure $\mu$, see \cite{royer06} and
\cite{khapha12}.

%This slightly extends the comparison Theorem 2.5 in \cite{royer06}, from which we borrow the main arguments.

\begin{Theorem}
\label{CompThm}
Let $\xi^{1},\xi^{2}\in{\bf L^{2}}(\Fc_{T})$ be two terminal conditions and let $F^{1},F^{2}:\Omega\times[0,T]\times\R\times\R^d \times \mathbf{L^2(\lambda)}\rightarrow \R$ be two generators satisfying  the assumptions 2.(i)-(iii) of Section 2.
Let $(Y^{1},Z^{1},L^{1},K^{1,-})$ $\in$ ${\bf S^2}\times{\bf L^2(W)}\times{\bf L^2(\tilde \mu)}\times{\bf K^2}$ satisfying
\beq
Y_t^1 & = & \xi^1 + \int_t^T F^1(s,Y_s^1,Z_s^1,L_s^1) ds - (K_T^{1,-}-K_t^{1,-})   \label{BSDEapp1}  \\
& &  -  \int_t^T  Z_s^1 dW_s  -  \int_t^T\int_A L_s^1(a)  \mu(ds,da)\;, \;\;\;\;\;   0 \leq t \leq T,  \; a.s. \nonumber
\enq
and $(Y^{2},Z^{2},L^{2},K^{2,+})$ $\in$ ${\bf S^2}\times{\bf L^2(W)}\times{\bf L^2(\tilde \mu)}\times{\bf K^2}$ satisfying
\beq
Y_t^2 & = & \xi^2 + \int_t^T F^2(s,Y_s^2,Z_s^2,L_s^2) ds + K_T^{2,+}-K_t^{2,+}   \label{BSDEapp2}  \\
& &  -  \int_t^T  Z_s^2 dW_s  -  \int_t^T\int_A L_s^2(a)  \mu(ds,da)\;, \;\;\;\;\;   0 \leq t \leq T,  \; a.s. \nonumber
\enq
If  $F^{1}(t,Y_{t}^{1},Z_{t}^{1},L_{t}^{1})$ $\leq$ $F^{2}(t,Y_{t}^{1},Z_{t}^{1},L_{t}^{1})$ $($resp.
$F^{1}(t,Y_{t}^{2},Z_{t}^{2},L_{t}^{2})$ $\leq$ $F^{2}(t,Y_{t}^{2},Z_{t}^{2},L_{t}^{2})$$)$,  $d\P\otimes dt$ a.e., and $\xi^{1} \leq \xi^{2}$ a.s., then
\beqs
Y_{t}^{1}\leq Y_{t}^{2}, \;\;\;\;\;\;   0 \leq t \leq T,  \; a.s.
\enqs
\end{Theorem}

\vspace{3mm}

We now state  a comparison theorem between a Skorohod solution and a Skorohod supersolution, both driven by the Brownian motion $W$ and the Poisson random measure $\mu$.
This slightly extends Theorem 5.2 in \cite{essaky}.

\begin{Theorem}
\label{CompThm2}
Let $\xi^{1},\xi^{2}\in{\bf L^{2}}(\Fc_{T})$ be two terminal conditions and let $F^{1},F^{2}:\Omega\times[0,T]\times\R\times\R^d \times \mathbf{L^2(\lambda)}\rightarrow \R$ be two generators satisfying assumptions 2.(i)-(iii) of Section 2.
Let $(Y^{1},Z^{1},L^{1},K^{1,-})$ $\in$ ${\bf S^2}\times{\bf L^2(W)}\times{\bf L^2(\tilde \mu)}\times{\bf K^2}$ satisfying
\beq
Y_t^1 & = & \xi^1 + \int_t^T F^1(s,Y_s^1,Z_s^1,L_s^1) ds - (K_T^{1,-}-K_t^{1,-})   \label{BSDEapp5}  \\
& &  -  \int_t^T  Z_s^1 dW_s  -  \int_t^T\int_A L_s^1(a)  \mu(ds,da)\;, \;\;\;\;\;   0 \leq t \leq T,  \; a.s. \nonumber
\enq
and
\beqs
Y_t^1 & \leq & U_t\;, \;\;\;\;\; 0 \leq t \leq T,\;a.s. \\
\int_0^T (U_{t^-} - Y_{t^-}^1) dK^{1,-}_t & = & 0\;, \;\;\;\;\;\;\, a.s.
\enqs
Furthermore, let $(Y^{2},Z^{2},L^{2},K^{2,+},K^{2,-})$ $\in$ ${\bf S^2}\times{\bf L^2(W)}\times{\bf L^2(\tilde \mu)}\times{\bf K^2}\times{\bf K^2}$ satisfying
\beq
Y_t^2 & = & \xi^2 + \int_t^T F^2(s,Y_s^2,Z_s^2,L_s^2) ds + K_T^{2,+}-K_t^{2,+} - (K_T^{2,-}-K_t^{2,-})   \label{BSDEapp6}  \\
& &  -  \int_t^T  Z_s^2 dW_s  -  \int_t^T\int_A L_s^2(a)  \mu(ds,da)\;, \;\;\;\;\;   0 \leq t \leq T,  \; a.s. \nonumber
\enq
and
\beqs
Y_t^2 & \leq & U_t\;, \;\;\;\;\; 0 \leq t \leq T,\;a.s. \\
\int_0^T (U_{t^-} - Y_{t^-}^2) dK^{2,-}_t & = & 0\;, \;\;\;\;\;\;\, a.s.
\enqs
If $\xi^{1} \leq \xi^{2}$ a.s. and $F^{1}(t,Y_{t}^{1},Z_{t}^{1},L_{t}^{1}) \leq F^{2}(t,Y_{t}^{1},Z_{t}^{1},L_{t}^{1})$, $d\P\otimes dt$ a.e., then
\beqs
Y_{t}^{1}\leq Y_{t}^{2}, \;\;\;\;\;\;   0 \leq t \leq T,  \; a.s.
\enqs
\end{Theorem}
{\bf Proof.} Consider the following penalized BSDEs:
\beqs
Y_t^{n,1} &=&   \xi^1 + \int_t^T F^1(s,Y_s^{n,1},Z_s^{n,1},L_s^{n,1}) ds - n \int_t^T (U_s - Y_s^{n,1})^- ds \\
& & -  \int_t^T  Z_s^{n,1} dW_s  -  \int_t^T\int_A  L_s^{n,1}(a)  \mu(ds,da)
\enqs
and
\beqs
Y_t^{n,2} &=&   \xi^2 + \int_t^T F^2(s,Y_s^{n,2},Z_s^{n,2},L_s^{n,2}) ds + K_T^{2,+}-K_t^{2,+} - n \int_t^T (U_s - Y_s^{n,2})^- ds \\
& & -  \int_t^T  Z_s^{n,2} dW_s  -  \int_t^T\int_A  L_s^{n,2}(a)  \mu(ds,da),
\enqs
for all $0 \leq t \leq T$, almost surely. By comparison Theorem \ref{CompThm} we get $Y_{t}^{n,1}\leq Y_{t}^{n,2}$, for all $n\in\N$. Recalling Remark \ref{R:HamEssaky}, we have that $Y_t^{n,1}$ converges to $Y_t^1$. It remains to prove the convergence of $Y_t^{n,2}$ towards $Y_t^2$.

Set $\tilde Y^{n,2}$ $:=$ $Y^{n,2} + K^{2,+}$, $\tilde U$ $:=$ $U + K^{2,+}$, $\tilde\xi^2$ $:=$ $\xi^2 + K_T^{2,+}$, and $\tilde F^2(t,y,z,\ell)$ $:=$ $F^2(t,y-K_t^{2,+},z,\ell)$, for all $0 \leq t \leq T$, $y\in\R$, $z\in\R^d$, $\ell\in{\bf L^2(\lambda)}$, almost surely. Then
\beqs
\tilde{Y}_t^{n,2} &=&   \tilde{\xi}^2 + \int_t^T \tilde{F}^2(s,\tilde{Y}_s^{n,2},Z_s^{n,2},L_s^{n,2}) ds - n \int_t^T (\tilde{U}_s - \tilde{Y}_s^{n,2})^- ds \\
& & -  \int_t^T  Z_s^{n,2} dW_s  -  \int_t^T\int_A  L_s^{n,2}(a)  \mu(ds,da),
\enqs
for all $0 \leq t \leq T$, almost surely. Note that $\tilde\xi^2$ verifies the square integrability condition and $\tilde F^2$ satisfies assumptions 2.(i)-(iii) of Section 2. Moreover, $\tilde U_T\in {\bf S^2}$ and $\tilde U_T$ $\geq$ $\tilde\xi^2$, almost surely. Now, again from Remark \ref{R:HamEssaky}, we have that $\tilde Y^{n,2}$ converges to $\tilde Y^2$ $=$ $Y^2 + K^{2,+}$, and hence $Y^{n,2}$ converges to $Y^2$.
\ep

\setcounter{equation}{0} \setcounter{Assumption}{0}
\setcounter{Theorem}{0} \setcounter{Proposition}{0}
\setcounter{Corollary}{0} \setcounter{Lemma}{0}
\setcounter{Definition}{0} \setcounter{Remark}{0}

\subsection{Monotonic limit theorem for BSDEs with jumps}

We state  a monotonic limit theorem for BSDEs driven by the Brownian motion $W$ and the Poisson random measure $\mu$. This extends the monotonic limit Theorem~3.1 in \cite{pengxu05} to the jump case.

\begin{Theorem}
\label{MonLimThm}
Let $(Y^m,Z^m,L^m,K^{m,+},K^{m,-})_m$ be a sequence in ${\bf S^2}\times{\bf L^2(W)}\times{\bf L^2(\tilde \mu)}\times{\bf K^2}\times{\bf K^2}$, with $K^{m,+}$ continuous, solution to:
\beq
Y_t^m & = & \xi + \int_t^T F(s,Y_s^m,Z_s^m,L_s^m) ds    +  K_T^{m,+} - K_t^{m,+} - (K_T^{m,-} - K_t^{m,-})   \label{BSDEMonLimThm}  \\
& &  -  \int_t^T  Z_s^m dW_s  -  \int_t^T\int_A   L_s^m(a)  \mu(ds,da), \;\;\;\;\; 0 \leq t \leq T,\; a.s. \nonumber
\enq
such that
\beq \label{boundunim}
\sup_{m\in\N} \Big( \big\| Y^m \big\|_{_{\bf S^2}} + \big\| Z^m \big\|_{_{\bf L^2(W)}} + \big\| L^m \big\|_{_{\bf L^2(\tilde \mu)}} + \big\| K^{m,+} \big\|_{_{\bf S^2}} + \big\| K^{m,-} \big\|_{_{\bf S^2}} \Big) & < &  \infty,
\enq
and $(Y^m)_m$ converges increasingly to  $Y$ $\in$ ${\bf S^2}$.  Suppose also that the sequence $(K^{m,-})_m$ satisfies:
\beq
\label{K-croissantApp}
K_t^{m,-} - K_s^{m,-} &\leq& K_t^{m+1,-} - K_s^{m+1,-}, \;\;\;\;\; 0 \leq s \leq t \leq T,\;a.s.
\enq
for all $m\in\N$. Then there exists $(Z,L,K^{+},K^{-})\in{\bf L^2(W)}\times{\bf L^2(\tilde \mu)}\times{\bf K^2}\times{\bf K^2}$ such that
\beq
Y_t & = & \xi + \int_t^T F(s,Y_s,Z_s,L_s) ds    +  K_T^+ - K_t^+ - (K_T^- - K_t^-) \label{limitBSDE} \\
& &  -  \int_t^T  Z_s dW_s  -  \int_t^T\int_A   L_s(a)  \mu(ds,da), \;\;\;\;\; 0 \leq t \leq T,\; a.s. \nonumber
\enq
Here $(Z,L)$ is the strong $($resp. weak$)$ limit of $(Z^m,L^m)_m$ in ${\bf L^p(W)}\times{\bf L^p(\tilde \mu)}$, with $p\in[1,2)$, $($resp. in ${\bf L^2(W)}\times{\bf L^2(\tilde \mu)}$$)$. Furthermore, $K^+_t$ is the weak limit of $(K_t^{m,+})_m$ in ${\bf L^2}(\Fc_t)$, and $(K_t^{m,-})_m$ converges strongly up to $K^-_t$ in ${\bf L^2}(\Fc_t)$, for any $0 \leq t \leq T$.
\end{Theorem}
{\bf Proof.}
{\em Step 1.} {\em Limit BSDE.} From the boundedness condition \reff{boundunim}  and the Hilbert structure of ${\bf L^2(W)}\times{\bf L^2(\tilde \mu)}\times{\bf L^2(0,T)}$,  there exists a subsequence, $(Z^{m_k},L^{m_k},F(\cdot,Y^{m_k},Z^{m_k},$ $L^{m_k}))_k$ which converges weakly to some $(Z,L,G)\in{\bf L^2(W)}\times{\bf L^2(\tilde \mu)}\times{\bf L^2(0,T)}$. Thus, for each stopping time $\tau$ $\leq$ $T$, the following weak convergences hold in ${\bf L^2}(\Fc_\tau)$ as $k\rightarrow\infty$:
\beqs
\int_0^\tau F(s,Y_s^{m_k},Z_s^{m_k},L_s^{m_k})ds &\rightharpoonup& \int_0^\tau G(s)ds, \\
\int_0^\tau Z_s^{m_k} dW_s &\rightharpoonup& \int_0^\tau Z_s dW_s,\\
\int_0^\tau\int_A L_s^{m_k}(a) \mu(ds,da) &\rightharpoonup& \int_0^\tau \int_A L_s(a) \mu(ds,da).
\enqs
From \reff{K-croissantApp}, there exists $K^- \in {\bf K^2}$, such that  $K^-_t$ is the strong limit of $(K^{m_k,-}_t)_k$ in ${\bf L^2}(\Fc_t)$ for all $0\leq t\leq T$.
In particular,  $K^{m_k,-}_\tau \rightharpoonup K^-_\tau$. Moreover, since
\beqs
K_\tau^{m_k,+} & = & Y_0^{m_k} - Y_\tau^{m_k} + K^{m_k,-}_\tau - \int_0^\tau F(s,Y_s^{m_k},Z_s^{m_k},L_s^{m_k}) ds \\
& &\qquad\qquad\;\,  +  \int_0^\tau  Z_s^{m_k} dW_s  +  \int_0^\tau\int_A   L_s^{m_k}(a)  \mu(ds,da).
\enqs
we also have the weak convergence in ${\bf L^2}(\Fc_\tau)$
\beqs
K_\tau^{m_k,+} \;\,\rightharpoonup\;\, K_\tau^+ & := & Y_0 - Y_\tau + K^-_\tau - \int_0^\tau G(s) ds \\
& &\qquad\quad\;\;  +  \int_0^\tau  Z_s dW_s  +  \int_0^\tau\int_A   L_s(a)  \mu(ds,da),
\enqs
as $k\rightarrow\infty$. Note that $\E[(K^+_T)^2] < \infty$ and for any two stopping times $0 \leq \sigma \leq \tau \leq T$, we have $K^+_\sigma \leq K^+_\tau$ since $K^{m,+}_\sigma \leq K^{m,+}_\tau$. From this it follows that $K^+$ is an increasing process. Observe now that we have obtained the following decomposition for $Y$:
\beq
Y_t &=& Y_0 - \int_0^\tau G(s) ds - K_t^+ + K_t^-  +  \int_0^t  Z_s dW_s  +  \int_0^t\int_A   L_s(a)  \mu(ds,da). \label{Ydecomp}
\enq
Since the processes $K^{m_k,+}$ and $K^{m_k,-}$ are predictable, we deduce that $K^+$ and $K^-$ are also predictable. Besides, by Lemmas 3.1 and 3.2 of \cite{pengxu05}, $K^+$, $K^-$ and $Y$ are c\`adl\`ag processes. Thus, in the above decomposition of
$Y$ in  \reff{Ydecomp}, the components $Z$ and $L$ are unique. As a matter of fact, the uniqueness of $Z$ follows by identifying the Brownian parts and finite variation parts. The uniqueness of $L$ is then obtained by identifying the predictable parts and by recalling that the jumps of $\mu$ are totally inaccessible. From the uniqueness of $(Z,L)$, it  follows that the whole sequence $(Z^m,L^m)_m$ converges weakly to $(Z,L)$ in ${\bf L^2(W)}\times{\bf L^2(\tilde \mu)}$.

\vspace{3mm}

\ni {\em Step 2.} {\em Properties of the process $K^+$.} We establish that the contribution of the jumps of $K^+$ is mainly concentrated within a finite number of intervals with sufficiently small total length. More precisely, we apply Lemma 2.3 in \cite{peng00} to $K^+$.
Consequently, as in Lemma 2.3 in \cite{peng00}, for any $\delta,\varepsilon>0$, there exists a finite number of pairs of stopping times $(\sigma_k,\tau_k)$, $k=0,\ldots,N$, with $0 < \sigma_k \leq \tau_k \leq T$, such that all the intervals $(\sigma_k,\tau_k]$ are disjoint and
\beq
\E\sum_{k=0}^N (\tau_k - \sigma_k) \geq T - \frac{\varepsilon}{2}, \qquad\qquad \E\sum_{k=0}^N \sum_{\sigma_k < t \leq \tau_k} |\Delta K^+_t|^2 \leq \frac{\varepsilon\delta}{3}. \label{Kjumps}
\enq
We should note that in \cite{peng00} the filtration is Brownian, therefore it is continuous, and hence each stopping time $\sigma_k$ can be approximated by a sequence of announceable stopping times. In our case the stopping times $\sigma_k$'s are constructed as the successive times of jumps of the predictable process $K^+$ with size bigger than some given positive level, therefore each $\sigma_k$ is a predictable stopping time and the approximation of $\sigma_k$ by announceable stopping times is again possible. We can thus argue exactly the same way as in Lemma 2.3 in \cite{peng00} to derive both estimates in \reff{Kjumps}.

\vspace{3mm}

\ni {\em Step 3.} {\em Strong convergence.}
%We begin by noting that, since $K^{m,+}$ is continuous,
%\beqs
%\Delta(Y_t^m-Y_t) &=& \Delta K^+_t - \Delta K^-_t + \Delta K^{m,-}_t + \int_A \big(L_t^m(a) - L_t(a)\big)\mu(\{t\},da).
%\enqs
By  applying It\^o's formula to $|Y_t^m - Y_t|^2$ on a subinterval $(\sigma,\tau]$, with $0 \leq \sigma \leq \tau \leq T$, two stopping times, and recalling that $K^{m,+}$ is continuous, we obtain:
\beq
\E\big|Y_\tau^m - Y_\tau\big|^2 & = & \E\big|Y_\sigma^m - Y_\sigma\big|^2 + \E\int_\sigma^\tau |Z_s^m - Z_s|^2 ds +  \E\int_\sigma^\tau \int_A |L_s^m(a) - L_s(a)|^2 \lambda(da)ds
\nonumber  \\
& & \;  + \; 2\E\int_\sigma^\tau (Y_s^m - Y_s) \big(G(s) \!-\! F(s,Y_s^m,Z_s^m,L_s^m)\big)ds   \nonumber \\
& & \; + \; \E\!\sum_{t\in(\sigma,\tau]}\! |\Delta K^+_t \!-\! \Delta K^-_t \!+\! \Delta K^{m,-}_t|^2 \nonumber \\
& & \;  + \; 2\E\int_{(\sigma,\tau]} (Y_{s^-}^m - Y_{s^-}) dK^+_s - 2\E\int_{(\sigma,\tau]} (Y_{s^-}^m - Y_{s^-}) dK^-_s \nonumber \\
& & \;  - \; 2\E\int_{(\sigma,\tau]} (Y_{s}^m - Y_{s}) dK^{m,+}_s + 2\E\int_{(\sigma,\tau]} (Y_{s^-}^m - Y_{s^-}) dK^{m,-}_s  \nonumber \\
& & \;  + \; 2\E\int_{(\sigma,\tau]} \int_A \big(Y_{s}^m - Y_{s})(L_s^m(a) - L_s(a))  \lambda(da)ds.  \label{Ito}
\enq
Now, let us write
\beqs
\int_{(\sigma,\tau]} (Y_{s^-}^m - Y_{s^-}) dK^+_s &=&
\int_{(\sigma,\tau]} \big(Y_{s^-}^m + \Delta K^{m,-}_s - Y_{s^-} + \Delta K^+_s - \Delta K^-_s  \big) dK^+_s\\
& &  -  \sum_{t\in(\sigma,\tau]} (\Delta K^+_t)^2 +  \sum_{t\in(\sigma,\tau]} \Delta K^+_t \Delta(K^-_s - K^{m,-}_s),
\enqs
and observe  that
\beqs
\int_{(\sigma,\tau]} (Y_{s^-}^m - Y_{s^-}) d(K^-_s - K_s^{m,-}) \leq 0, & \mbox{ and } &
\int_{(\sigma,\tau]} (Y_{s}^m - Y_{s}) dK^{m,+}_s \leq 0.
\enqs
%while
%\begin{align*}
%2\E\int_{(\sigma,\tau]} \int_A \big(|Y_{s^-}^m - Y_{s^-} + L_s^m(a) &- L_s(a)|^2 - |L_s^m(a) - L_s(a)|^2\big) \mu(ds,da) \\
%&=\E\int_\sigma^\tau \int_A |L_s^m(a) - L_s(a)|^2 \lambda(da)ds \\
%&\;+ 2\E\int_\sigma^\tau \int_A (Y_{s^-}^m - Y_{s^-})(L_s^m(a) - L_s(a)) \lambda(da)ds.
%\end{align*}
Therefore, by using the inequality $2ab$ $\geq$ $-2b^2$ $-$ $a^2/2$, we obtain from \reff{Ito}
\beqs
 & &\E\int_\sigma^\tau |Z_s^m - Z_s|^2 ds + \frac{1}{2}\E\int_\sigma^\tau \int_A |L_s^m(a) - L_s(a)|^2 \lambda(da)ds \label{Ito2}  \\
& \leq &  \E\big|Y_\tau^m - Y_\tau\big|^2 + 2\lambda(A)\E\int_\sigma^\tau\big|Y_s^m - Y_s\big|^2ds \nonumber \\
& & \; + \;  2\E\int_\sigma^\tau \big|Y_s^m - Y_s\big| \big|G(s) - F(s,Y_s^m,Z_s^m,L_s^m)\big|ds \nonumber \\
& &\; - 2\E\int_{(\sigma,\tau]} \big(Y_{s^-}^m + \Delta K^{m,-}_s - Y_{s^-} + \Delta K^+_s - \Delta K^-_s  \big) dK^+_s + 2 \E\sum_{t\in(\sigma,\tau]} |\Delta K^+_t|^2\nonumber\\
& &\;  - 2 \E\sum_{t\in(\sigma,\tau]} \Delta K^+_t \Delta(K^-_s - K^{m,-}_s) - \E\sum_{t\in(\sigma,\tau]} |\Delta K^+_t - \Delta K^-_t + \Delta K^{m,-}_t|^2, \nonumber \\
& \leq &  \E\big|Y_\tau^m - Y_\tau\big|^2 + 2\lambda(A)\E\int_\sigma^\tau\big|Y_s^m - Y_s\big|^2ds \nonumber \\
& &  \;+ 2\E\int_\sigma^\tau \big|Y_s^m - Y_s\big| \big|G(s) - F(s,Y_s^m,Z_s^m,L_s^m)\big|ds \nonumber \\
& & \;  - 2\E\int_{(\sigma,\tau]} \big(Y_{s^-}^m + \Delta K^{m,-}_s - Y_{s^-} + \Delta K^+_s - \Delta K^-_s  \big) dK^+_s + \E\sum_{t\in(\sigma,\tau]} |\Delta K^+_t|^2. \nonumber
\enqs
by using  the inequality $2a^2 - 2ab -(a-b)^2 \leq a^2$.
We know that the first two terms on the right-hand side of \reff{Ito2} converge to zero as $m\rightarrow\infty$. The third term also tends to zero since $(G(\cdot)$ $-$ $F(\cdot,Y^m,Z^m,L^m))_m$ is bounded in ${\bf L^2(0,T)}$, and so by
Cauchy-Schwarz inequality
\beqs
\E\int_0^T \big|Y_s^m - Y_s\big| \big|G(s) - F(s,Y_s^m,Z_s^m,L_s^m)\big|ds  &  \rightarrow  &  0, \qquad \mbox{ as } m\rightarrow\infty.
\enqs
For the fourth term, since $K^{m,-}$ is predictable, the predictable projection of $Y^m$ is $^pY_t^m$ $=$ $Y_{t^-}^m + \Delta K_t^{m,-}$. Similarly, from \reff{Ydecomp} and since $K^+$ and $K^-$ are predictable processes, we see that $^pY_t$ $=$ $Y_{t^-} - \Delta K^+_t + \Delta K^-_t$. By the dominated convergence theorem, we obtain
\beqs
\lim_{m\rightarrow\infty} \E\int_{(\sigma,\tau]} \big(Y_{s^-}^m + \Delta K^{m,-}_s - Y_{s^-} + \Delta K^+_s - \Delta K^-_s  \big) dK^+_s &=&  0.
\enqs
For the last term in \reff{Ito2}, we exploit the results in \reff{Kjumps}, regarding the contribution of the jumps of $K^+$. More precisely, we apply estimate \reff{Ito2} for each $\sigma=\sigma_k$ and $\tau=\tau_k$, with $\sigma_k,\tau_k$ defined in Step 2, and then take the sum over $k=0,\ldots,N$. It follows that
\beqs
& &\sum_{k=0}^N\E\int_{\sigma_k}^{\tau_k} |Z_s^m - Z_s|^2 ds + \frac{1}{2}\sum_{k=0}^N\E\int_{\sigma_k}^{\tau_k} \int_A |L_s^m(a) - L_s(a)|^2 \lambda(da)ds\\
& \leq & \sum_{k=0}^N\E\big|Y_{\tau_k}^m - Y_{\tau_k}\big|^2 + 2\lambda(A)\E\int_0^T\big|Y_s^m - Y_s\big|^2ds \\
& &\;  + 2\E\int_0^T \big|Y_s^m - Y_s\big| \big|G(s) - F(s,Y_s^m,Z_s^m,L_s^m)\big|ds + \sum_{k=0}^N\E\sum_{t\in(\sigma_k,\tau_k]} |\Delta K^+_t|^2 \\
& & \;  - 2\sum_{k=0}^N\E\int_{(\sigma_k,\tau_k]} \big(Y_{s^-}^m + \Delta K^{m,-}_s - Y_{s^-} + \Delta K^+_s - \Delta K^-_s  \big) dK^+_s.
\enqs
From the above convergence results, we deduce that
\beqs
 & &\limsup_{m\rightarrow\infty} \bigg(\sum_{k=0}^N\E\int_{\sigma_k}^{\tau_k} |Z_s^m - Z_s|^2 ds + \frac{1}{2}\sum_{k=0}^N\E\int_{\sigma_k}^{\tau_k} \int_A |L_s^m(a) - L_s(a)|^2 \lambda(da)ds\bigg) \\
&\leq&  \sum_{k=0}^N\E\sum_{t\in(\sigma_k,\tau_k]} |\Delta K^+_t|^2 \leq \frac{\varepsilon\delta}{3}.
\enqs
Therefore, following the same steps as in the proof of Theorem 2.1 in \cite{peng00}, we deduce that the sequences $(Z^m)_m$ and $(L^m)_m$ converge in measure, respectively, to $Z$ and $L$. Since they are bounded, respectively, in ${\bf L^2(W)}$ and ${\bf L^2(\tilde \mu)}$, they are uniformly integrable in ${\bf L^p(W)}$ and ${\bf L^p(\tilde \mu)}$, for any $p\in[1,2)$. Thus, $(Z^m)_m$ and $(L^m)_m$ converge strongly to $Z$ and $L$ in ${\bf L^p(W)}$ and ${\bf L^p(\tilde \mu)}$, respectively.

By the Lipschitz condition on $F$, we also have the strong convergence in ${\bf L^p(0,T)}$ of $(F(\cdot,Y^m,Z^m,L^m))_m$ to $F(\cdot,Y,Z,L)$. Since $G(\cdot)$ is the weak limit of $(F(\cdot,Y^m,Z^m,L^m))_m$ in ${\bf L^2(0,T)}$, we deduce that $G(\cdot)$ $=$ $F(\cdot,Y,Z,L)$. Therefore we obtain that $(Y,Z,L,K^+,K^-)$ satisfies the BSDE \reff{limitBSDE}.
\ep

\end{document}